\documentclass{article}
\usepackage{amsfonts}
\usepackage{amsmath}

\newtheorem{theorem}{Theorem}

\newtheorem{corollary}[theorem]{Corollary}

\newtheorem{lemma}[theorem]{Lemma}

\begin{document}

\title{Some functionals for random walks and critical branching processes in
extreme random environment}
\author{Congzao Dong\thanks{%
Xidian University, 266 Xinglong Section of Xifeng Road, Xi'an, Shaanxi,
710126, China Email: czdong@xidian.edu.cn}, \ Elena Dyakonova\thanks{%
Steklov Mathematical Institute of Russian Academy of Sciences, 8 Gubkina
St., Moscow 119991 Russia Email: elena@mi-ras.ru}, and Vladimir Vatutin
\thanks{%
Steklov Mathematical Institute of Russian Academy of Sciences, 8 Gubkina
St., Moscow 119991 Russia Email: vatutin@mi-ras.ru}}
\date{}
\maketitle

\begin{abstract}
Let $\left\{ S_{n},n\geq 0\right\} $ be a random walk whose increment
distribution belongs without centering to the domain of attraction of an $%
\alpha $-stable law, i.e., there are some scaling constants $a_{n}$ such
that the sequence $S_{n}/a_{n},n=1,2,...,$ weakly converges, as $%
n\rightarrow \infty $ to a random variable having an $\alpha $-stable
distribution. Let $S_{0}=0,$%
\begin{equation*}
L_{n}:=\min \left( S_{1},...,S_{n}\right) ,\tau _{n}:=\min \left\{ 0\leq
k\leq n:S_{k}=\min (0,L_{n})\right\} .
\end{equation*}%
Assuming that $S_{n}\leq h(n),$ where $h(n)$ is $o(a_{n})$ and $%
\lim_{n\rightarrow \infty }h(n)\in \lbrack -\infty ,+\infty ]$ exists we
prove several limit theorems describing the asymptotic behavior of the
functionals
\begin{equation*}
\mathbf{E}\left[ e^{S_{\tau _{n}}};S_{n}\leq h(n)\right]
\end{equation*}%
as $n\rightarrow \infty $. The obtained results are applied for studying the
survival probability of a critical branching process evolving in an
extremely unfavorable random environment.


\end{abstract}

\section{Introduction and main results}

We consider some functionals specified on the trajectories of a random walk
\begin{equation*}
S_{0}=0,\quad S_{n}=X_{1}+...+X_{n},\ n\geq 1,
\end{equation*}%
with i.i.d. increments $X_{i},i=1,2,...$. To describe the conditions we
impose on the increments, let
\begin{equation*}
\mathcal{A}:=\{\alpha \in (0,2)\backslash \{1\},\,|\beta |<1\}\cup \{\alpha
=1,\beta =0\}\cup \{\alpha =2,\beta =0\}
\end{equation*}%
be a subset in $\mathbb{R}^{2}.$ For $(\alpha ,\beta )\in \mathcal{A}$ and a
random variable $X$ we write $X\in \mathcal{D}\left( \alpha ,\beta \right) $
if the distribution of $X$ belongs to the domain of attraction of a stable
law with density $g_{\alpha ,\beta }(x),x\in (-\infty ,+\infty ),$ and the
characteristic function%
\begin{equation*}
G_{\alpha ,\beta }(w)=\int_{-\infty }^{+\infty }e^{iwx}g_{\alpha .\beta
}(x)\,dx=\exp \left\{ -c|w|^{\,\alpha }\left( 1-i\beta \frac{w}{|w|}\tan
\frac{\pi \alpha }{2}\right) \right\} ,\ c>0,
\end{equation*}%
and, in addition, $\mathbf{E}X=0$ if this moment exists. This implies, in
particular, that there is an increasing sequence of positive numbers
\begin{equation*}
a_{n}\ =\ n^{1/\alpha }\ell (n)
\end{equation*}%
with a slowly varying sequence $\ell (1),\ell (2),\ldots ,$ such that, as $%
n\rightarrow \infty $
\begin{equation*}
\left\{ \frac{S_{\left[ nt\right] }}{a_{n}},t\geq 0\right\} \Longrightarrow
\mathcal{Y}=\left\{ Y_{t},t\geq 0\right\} ,
\end{equation*}%
where%
\begin{equation*}
\mathbf{E}e^{iwY_{t}}=G_{\alpha ,\beta }(wt^{1/\alpha }),t\geq 0,
\end{equation*}%
and the symbol $\Longrightarrow $ stands for the weak convergence in the
space $D[0,\infty )$ of c\`{a}dl\`{a}g functions endowed with Skorokhod
topology. Observe that if $X_{n}\overset{d}{=}X\in \mathcal{D}\left( \alpha
,\beta \right) $ for all $n\in \mathbb{N}:\mathbf{=}\left\{ 1,2,...\right\} $
then%
\begin{equation*}
\lim_{n\rightarrow \infty }\mathbf{P}\left( S_{n}>0\right) =\rho =\mathbf{P}%
\left( Y_{1}>0\right) \in (0,1).
\end{equation*}

We now list our main restrictions on the properties of the random walk.

\textbf{Condition A1.} \textit{The random variables }$X_{n},n\in \mathbb{N},$%
\textit{\ are independent copies of a random variable }$X\in \mathcal{D}%
\left( \alpha ,\beta \right) $\textit{.\ Besides, the distribution of }$X$%
\textit{\ is non-lattice.}\emph{\ }

Some our statements need a stronger assumption.

\textbf{Condition A2. }\textit{The law of }$X$\textit{\ under }$\mathbf{P}$%
\textit{\ is absolutely continuous with respect to the Lebesgue measure on }$%
\mathbb{R}$\textit{, and there exists }$n\in \mathbb{N}$\textit{\ such that
the density }$f_{n}(x):=\mathbf{P}(S_{n}\in dx)/dx$\textit{\ of }$S_{n}$%
\textit{\ is essentially bounded (therefore, }$f_{n}(x)\in L^{\infty }$%
\textit{).}

To formulate the main results of the paper we denote%
\begin{eqnarray}
M_{n} &:=&\max \left( S_{1},...,S_{n}\right) ,\quad  \notag \\
&&  \notag \\
L_{n} &:=&\min \left( S_{1},...,S_{n}\right) ,\tau _{n}:=\min \left\{ 0\leq
k\leq n:S_{k}=\min (0,L_{n})\right\} ,  \notag \\
b_{n} &:=&\frac{1}{a_{n}n}=\frac{1}{\ n^{1/\alpha +1}\ell (n)},  \label{Defb}
\end{eqnarray}%
and introduce renewal functions%
\begin{eqnarray*}
U(x) &:=&I\left\{ x\geq 0\right\} +\sum_{n=1}^{\infty }\mathbf{P}\left(
S_{n}\geq -x,M_{n}<0\right) \\
&=&I\left\{ x\geq 0\right\} +\sum_{n=1}^{\infty }\mathbf{P}\left( S_{n}\geq
-x,\tau _{n}=n\right) ,\ x\in \mathbb{R},
\end{eqnarray*}%
\begin{equation*}
V_{0}(x)\ :=\ I\left\{ x\leq 0\right\} +\sum_{n=1}^{\infty }\mathbf{P}\left(
S_{n}\leq -x,L_{n}\geq 0\right) ,\ x\in \mathbb{R},
\end{equation*}%
and

\begin{equation}
V(x)\ :=\ I\left\{ x<0\right\} +\sum_{n=1}^{\infty }\mathbf{P}\left(
S_{n}<-x,L_{n}\geq 0\right) ,\ x\in \mathbb{R}.  \label{RenUV}
\end{equation}%
Observe that
\begin{equation}
V_{0}(0)=1+\sum_{n=1}^{\infty }\mathbf{P}\left( S_{n}=0,L_{n}\geq 0\right) =%
\frac{1}{1-\zeta },  \label{V_zero}
\end{equation}%
where%
\begin{eqnarray*}
\zeta &=&\mathbf{P}\left( S_{1}=0\right) +\sum_{n=2}^{\infty }\mathbf{P}%
\left( S_{1}>0,...,S_{n-1}>0,S_{n}=0\right) \\
&=&\mathbf{P}\left( S_{1}=0\right) +\sum_{n=2}^{\infty }\mathbf{P}\left(
S_{1}<0,...,S_{n-1}<0,S_{n}=0\right) \in (0,1).
\end{eqnarray*}%
Here to justify the transition from the first line to the second line it is necessary to use the fact that
\begin{equation*}
\left\{ S_{n}-S_{n-k},k=0,1,...,n\right\} \overset{d}{=}\left\{
S_{k},k=0,1,...,n\right\} .
\end{equation*}%
One can check also that
\begin{equation*}
U(0+)=V(-0)=1,
\end{equation*}%
and $V(x)=V_{0}(x)$ if Condition A2 is valid.

\begin{theorem}
\label{T_minimPosit} Let Condition A1 be valid. If $\varphi (n),n\in \mathbb{%
N}\mathbf{,}$ is a positive deterministic function such that $\varphi
(n)\rightarrow +\infty $ as $n\rightarrow \infty $ and $\varphi
(n)=o(a_{n}), $ then, \ for any $\theta >0$%
\begin{eqnarray*}
\mathbf{E}\left[ e^{\theta S_{\tau _{n}}};S_{n}\leq \varphi (n)\right] &\sim
&\theta \mathbf{P}(S_{n}\leq \varphi (n),L_{n}\geq 0)\int_{0}^{\infty
}e^{-\theta z}U(z)dz \\
&\sim &\theta g_{\alpha ,\beta }(0)b_{n}\int_{0}^{\varphi
(n)}V(-w)dw\int_{0}^{\infty }e^{-\theta z}U(z)dz
\end{eqnarray*}%
as $n\rightarrow \infty $.
\end{theorem}

We now formulate three statements for the case when $S_{n}$ $\rightarrow
-\infty $ as $n\rightarrow \infty .$

\begin{theorem}
\label{T_lowend} Let Condition A1 be valid. If $\psi (n),n\in \mathbb{N},$
is a negative deterministic function such that $\psi (n)\rightarrow -\infty $
as $n\rightarrow \infty $ and $\psi (n)=o(a_{n}),$ then, for any $\theta >0$
and $x\leq 0$%
\begin{equation}
\mathbf{E}_{x}\left[ e^{\theta S_{n}};S_{n}\leq \psi (n),M_{n}<0\right] \sim
g_{\alpha ,-\beta }(0)b_{n}V(x)U\left( -\psi (n)\right) \theta
^{-1}e^{\theta \psi (n)}  \label{MaximEnd}
\end{equation}%
as $n\rightarrow \infty $.
\end{theorem}

Using the duality principle for random walks and setting $x=0$ in (\ref%
{MaximEnd}) we immediately get the following result:

\begin{corollary}
\label{C_VatVat} If Condition A1 is valid and $\psi (n),n\in \mathbb{N}%
\mathbf{,}$ is a negative deterministic function such that $\psi
(n)\rightarrow -\infty $ as $n\rightarrow \infty $ and $\psi (n)=o(a_{n}),$
then, for any $\theta >0$
\begin{equation}
\mathbf{E}\left[ e^{\theta S_{n}};S_{n}\leq \psi (n),\tau _{n}=n\right] \sim
g_{\alpha ,-\beta }(0)b_{n}U\left( -\psi (n)\right) \theta ^{-1}e^{\theta
\psi (n)}  \label{Min_end}
\end{equation}%
as $n\rightarrow \infty $.
\end{corollary}

The next statement is a natural complement to Corollary \ref{C_VatVat}.

\begin{theorem}
\label{T_minim} Let Condition A1 be valid. If $\psi (n),n\in \mathbb{N},$ is
a negative deterministic function such that $\psi (n)\rightarrow -\infty $
as $n\rightarrow \infty $ and $\psi (n)=o(a_{n}),$ then, \ for any $\theta
>0 $%
\begin{equation*}
\mathbf{E}\left[ e^{\theta S_{\tau _{n}}};S_{n}\leq \psi (n)\right] \sim
g_{\alpha ,-\beta }(0)b_{n}U\left( -\psi (n)\right) e^{\theta \psi
(n)}\int_{0}^{\infty }e^{-\theta z}V_{0}(-z)dz.
\end{equation*}
\end{theorem}

We now consider the case $S_{n}\leq K$ for a fixed $K$.

\begin{theorem}
\label{T_interm} Let Conditions A1 and A2 be valid. Then, for any fixed $K$
and $\theta >0$
\begin{eqnarray*}
\lim_{n\rightarrow \infty }\frac{\mathbf{E}\left[ e^{\theta S_{\tau
_{n}}};S_{n}\leq K\right] }{b_{n}} &=&g_{\alpha ,\beta }(0)\int_{-\infty
}^{K\wedge 0}e^{\theta x}U(-dx)\int_{0}^{K-x}V(-w)dw \\
&&+g_{\alpha ,-\beta }(0)\int_{-\infty }^{K\wedge 0}e^{\theta
x}U(-x)V_{0}(K-x)dx.
\end{eqnarray*}
\end{theorem}

We note that the problems related to those we consider here were
investigated by several authors.

Thus, Hirano \cite{Hir98} analyzed the asymptotic behavior, as $n\rightarrow
\infty $ of the functional%
\begin{equation*}
\mathbf{E}_{x}\left[ e^{\theta S_{n}};M_{n}\leq K\right] =e^{\theta x}%
\mathbf{E}\left[ e^{\theta S_{n}};M_{n}\leq K-x\right]
\end{equation*}%
assuming that $\mathbf{E}X^{2}<\infty $ and the pair $(K,x)$ is fixed.
Afanasyev et al \cite{ABGV2011} generalized Hirano's results to the case $%
X\in \mathcal{D}(\alpha ,\beta )$ by showing, as $n\rightarrow \infty $
that, for $\theta >0$ and fixed $x\leq 0$
\begin{equation}
\mathbf{E}_{x}\left[ e^{\theta S_{n}};M_{n}<0\right] \sim g_{\alpha ,\beta
}(0)b_{n}V(x)\int_{0}^{\infty }e^{-\theta z}U(z)dz,  \label{MaxSmall}
\end{equation}%
while for fixed $x\geq 0$%
\begin{equation*}
\mathbf{E}_{x}\left[ e^{-\theta S_{n}};L_{n}\geq 0\right] \sim g_{\alpha
,\beta }(0)b_{n}U(x)\int_{0}^{\infty }e^{-\theta z}V(-z)dz.
\end{equation*}%
It follows from (\ref{MaxSmall}) and the continuity theorem for Laplace
transforms that, for fixed $y\leq 0$
\begin{equation*}
\mathbf{E}\left[ e^{\theta S_{n}}I\left\{ S_{n}<y\right\} ;M_{n}<0\right]
\sim g_{\alpha ,\beta }(0)b_{n}\int_{-y}^{\infty }e^{-\theta z}U(z)dz
\end{equation*}%
as $n\rightarrow \infty $ which gives a hint for the form of the asymptotic
behavior of the left-hand side in (\ref{Min_end}).

The structure of the remaining sections of the paper looks as follows. In
Section \ref{Sec2} we formulate a number of known results for random walks
conditioned to stay nonnegative or negative. In Section \ref{Sec3} we prove
Theorem \ref{T_minimPosit}. Section \ref{Sec4} is devoted to the proof of
Theorem \ref{T_lowend}. Section \ref{Sec5} contains the proof of Theorem \ref%
{T_minim}. The proof of Theorem \ref{T_interm} is given in Section \ref{Sec6}%
. In Section \ref{Sec7} we introduce measures $\mathbf{P}_{x}^{+}$ and $%
\mathbf{P}_{x}^{-}$ generated, respectively, by random walks conditioned to
stay nonegative or negative and use Theorem \ref{T_interm} to investigate
the survival probability of a critical branching process evolving in extreme
random environment.

In what follows we denote by $C,C_{1},C_{2},...,$ some positive constants
that may be different in different formulas or even within one and the same
formula.

\section{Auxilary results\label{Sec2}}

In the sequel we need to consider random walks that start from any point $%
x\in \mathbb{R}$ or from an initial distribution $\nu $. In such cases we
write for probabilities as usual $\mathbf{P}_{x}\left( \cdot \right) $ or $%
\mathbf{P}_{\nu }\left( \cdot \right) $. We also write $\mathbf{P}$ instead
of $\mathbf{P}_{0}$.

We now formulate a number of statements that show importance of the
functions $U,V_{0}$ and $V$.

We recall that a positive sequence $\left\{ c_{n},n\in \mathbb{N}\right\} $
--- or a real function $c(x),x\geq x_{0},$ --- is said to be regularly
varying at infinity with index $\gamma \in \mathbb{R}$ , denoted $c_{n}\in
R_{\gamma }$ or $c(x)\in R_{\gamma }$ if $c_{n}\sim n^{\gamma }l(n)$ $%
(c(x)\sim x^{\gamma }l(x))$, where $l(x)$ is a slowly varying function, i.e.
a positive real function with the property that $l(tx)/l(x)\rightarrow 1$ as
$x\rightarrow \infty $ for all fixed $t>0$.

It is known (see, for instance, \cite{Rog1971}, \cite{Sin57}) that if
Condition A1 is valid then
\begin{equation}
\mathbf{P}\left( M_{n}<0\right) \in R_{-\rho },\quad U(x)\in R_{\alpha \rho
},  \label{AsymU}
\end{equation}%
and%
\begin{equation}
\mathbf{P}\left( L_{n}\geq 0\right) \in R_{-(1-\rho )},\quad V(-x)\in
R_{\alpha (1-\rho )}.  \label{AsymV}
\end{equation}

Some basic inequalities used below in our proofs are contained in the
following lemma.

\begin{lemma}
\label{L_double}(see \cite[Proposition 2.3]{ABGV2011} and \cite[Lemma 2]%
{VD2022}) If Condition A1 is valid, then there is a number $C>0$ such that,
for all $n$ and $x,y\geq 0$,
\begin{equation}
\mathbf{P}_{x}\left( 0\leq S_{n}<y,L_{n}\geq 0\right) \ \leq \
C\,b_{n}\,U(x)\int_{0}^{y}V(-z)dz\ ,  \label{Rough1}
\end{equation}%
and, for $x,z\leq 0$,
\begin{equation}
\mathbf{P}_{x}\left( z\leq S_{n}<z+1,M_{n}<0\right) \ \leq \
C\,b_{n}\,V(x)U(-z),  \label{LocalMaxim}
\end{equation}%
and, therefore,
\begin{equation*}
\mathbf{P}_{x}\left( z\leq S_{n}<0,M_{n}<0\right) \ \leq \
C\,b_{n}\,V(x)\int_{z}^{0}U(-w)\ dw.
\end{equation*}
\end{lemma}

We need some equivalence relations established by Doney \cite{Don12} and
rewritten in our notation.

\begin{lemma}
\label{L_DonLocal} (Proposition 18 in \cite{Don12}) Suppose that $X\in
\mathcal{D}\left( \alpha ,\beta \right) $ and the distribution of $X$ is
non-lattice. Then, for each fixed $\Delta >0$%
\begin{equation}
\mathbf{P}_{z}\left( S_{n}\in \lbrack y,y+\Delta ),L_{n}\geq 0\right) \sim
g_{\alpha ,\beta }(0)b_{n}U(z)\int_{y}^{y+\Delta }V(-w)dw  \label{DDon0}
\end{equation}%
and%
\begin{equation}
\mathbf{P}_{-z}\left( S_{n}\in \lbrack -y,-y+\Delta ),M_{n}<0\right) \sim
g_{\alpha ,-\beta }(0)b_{n}V(-z)\int_{y-\Delta }^{y}U(w)dw  \label{DDon1}
\end{equation}
uniformly in $\max (z,y)\in \lbrack 0,\delta _{n}a_{n}]$, where $\delta
_{n}\rightarrow 0$ as $n\rightarrow \infty $.
\end{lemma}

Integration (\ref{DDon0}) over $y\in (0,x)$ leads to the following important
conclusion for $z=0$ (see also, Theorem 4 in \cite{VW2010}).

\begin{corollary}
\label{C_IntegVW} Under the conditions of Lemma \ref{L_DonLocal}
\begin{equation*}
\mathbf{P(}S_{n}\leq x,L_{n}\geq 0)\sim g_{\alpha ,\beta
}(0)b_{n}\int_{0}^{x}V(-w)dw
\end{equation*}%
uniformly in \thinspace $x\in (0,\delta _{n}a_{n}],$ where $\delta
_{n}\rightarrow 0$ as $n\rightarrow \infty $.
\end{corollary}

Combining (\ref{AsymU}) and (\ref{DDon1}) we arrive to the following
statement.

\begin{corollary}
\label{C_Doney_regular} Suppose that $X\in \mathcal{D}\left( \alpha ,\beta
\right) $ and the distribution of $X$ is non-lattice. Then, for each fixed $%
\Delta >0$
\begin{equation*}
\mathbf{P}_{-x}\left( S_{n}\in \lbrack -y,-y+\Delta ),M_{n}<0\right) \sim
g_{\alpha ,-\beta }(0)b_{n}V(-x)U(y)\Delta
\end{equation*}
uniformly in $x\in \lbrack 0,\delta _{n}a_{n}]$ and $y\in \lbrack
T_{n},\delta _{n}a_{n}]$, where $T_{n}\rightarrow \infty $ and $\delta
_{n}\rightarrow 0$ as $n\rightarrow \infty $ in such a way that $%
T_{n}<\delta _{n}a_{n}$.
\end{corollary}

We also need the following simple observation, which we will refer to
several times in the sequel.

\begin{lemma}
\label{L_inequal} Let $g(x),x\geq 0,$ be a positive nondecreasing function
such that $g(2x_{0})\geq 1$ for some $x_{0}>0$. Then, for any $y\geq 0$%
\begin{equation*}
g(x_{0}+y)\leq g(2x_{0})(1+g(2y)).
\end{equation*}%
If, in addition, $g(x)\rightarrow \infty $ as $x\rightarrow \infty $ and $%
g(x)\in R_{\gamma }$ for some $\gamma \in \mathbb{R},$ then there is a
constant $C\in (0,\infty )$ such that%
\begin{equation}
g(x+y)\leq Cg(x)(1+g(2y))  \label{DoubleInequ}
\end{equation}%
for all $y\geq 0$ and all sufficiently large $x$.
\end{lemma}

\section{Proof of Theorem \ref{T_minimPosit}}

\label{Sec3}

We fix $J<n$ and write%
\begin{equation}
\mathbf{E}\left[ e^{\theta S_{\tau _{n}}};S_{n}\leq \varphi (n)\right]
=\sum_{j=0}^{J}\mathbf{E}\left[ e^{\theta S_{j}};S_{n}\leq \varphi (n),\tau
_{n}=j\right] +R(J+1,n),  \label{Decomp1}
\end{equation}%
where%
\begin{equation*}
R(J+1,n):=\sum_{j=J+1}^{n}\mathbf{E}\left[ e^{\theta S_{j}};S_{n}\leq
\varphi (n),\tau _{n}=j\right] .
\end{equation*}%
It follows from (\ref{LocalMaxim}) and (\ref{Rough1}) that, for all positive
integers $n$ and $k$
\begin{eqnarray}
\mathbf{P}\left( S_{n}\in \lbrack -k,-k+1),M_{n}<0\right) &=&\mathbf{P}%
\left( S_{n}\in \lbrack -k,-k+1),\tau _{n}=n\right)  \notag \\
&\leq &C\,b_{n}\,U(k)  \label{Rough00}
\end{eqnarray}%
and
\begin{equation}
\mathbf{P}\left( 0\leq S_{n}<y,L_{n}\geq 0\right) \ \leq \
C\,b_{n}\,\int_{0}^{y}V(-z)dz.  \label{Rough2}
\end{equation}%
Let $\left\{ S_{n}^{\prime },n\geq 0\right\} $ be an independent
probabilistic copy of the random walk\newline
 $\left\{ S_{n},n\geq 0\right\} $ and $%
L_{n}^{\prime }=\min \left\{ S_{1}^{\prime },...,S_{n}^{\prime }\right\} $.
Using (\ref{Rough00}) and (\ref{Rough2}) we obtain
\begin{eqnarray}
&&\mathbf{E}\left[ e^{\theta S_{j}};S_{n}\leq \varphi (n),\tau _{n}=j\right]
\notag \\
&=&\mathbf{E}\left[ e^{\theta S_{j}}\mathbf{P}\left( S_{n-j}^{\prime }\leq
\varphi (n)-S_{j},L_{n-j}^{\prime }\geq 0|S_{j}\right) ;\tau _{j}=j\right]
\notag \\
&\leq &\sum_{k=1}^{\infty }e^{\theta \left( -k+1\right) }\mathbf{P}\left(
S_{j}\in \lbrack -k,-k+1),\tau _{j}=j\right) \mathbf{P}\left( S_{n-j}\leq
\varphi (n)+k,L_{n-j}\geq 0\right)  \notag \\
&\leq &Cb_{j}\sum_{k=1}^{\infty }e^{\theta \left( -k+1\right) }U(k)\mathbf{P}%
\left( S_{n-j}\leq \varphi (n)+k,L_{n-j}\geq 0\right)  \label{RightTail0} \\
&\leq &Cb_{j}b_{n-j}\sum_{k=1}^{\infty }e^{-\theta k}U(k)\int_{0}^{\varphi
(n)+k}V(-z)dz.\   \label{RightTail}
\end{eqnarray}%
Note that in view of (\ref{Defb}) and properties of regularly varying
functions there exists a constant $C\in (0,\infty )$ such that $b_{j}\leq
Cb_{n}$ for all $n$ and $j\in \lbrack n/2,n]$.

Since $U(x)$ and $V(-x)$ are $O(|x|+1)$ as renewal functions, the estimate%
\begin{equation*}
\sum_{k=1}^{\infty }e^{-\theta k}U(k)(1+V(-2k))<\mathbf{\infty }
\end{equation*}%
is valid for all $\theta >0$. Hence we conclude by (\ref{DoubleInequ}) that
\begin{eqnarray}
&&\sum_{j=\left[ n/2\right] +1}^{n}\mathbf{E}\left[ e^{\theta
S_{j}};S_{n}\leq \varphi (n),\tau _{n}=j\right]  \notag \\
&&\qquad \leq C\sum_{k=1}^{\infty }e^{-\theta k}U(k)\sum_{j=\left[ n/2\right]
+1}^{n}b_{j}\mathbf{P}\left( S_{n-j}\leq \varphi (n)+k,L_{n-j}\geq 0\right)
\notag \\
&&\qquad \leq C_{1}b_{n}\sum_{k=1}^{\infty }e^{-\theta k}U(k)\sum_{j=\left[
n/2\right] +1}^{n}\mathbf{P}\left( S_{n-j}\leq \varphi (n)+k,L_{n-j}\geq
0\right)  \notag \\
&&\qquad \leq C_{1}b_{n}\sum_{k=1}^{\infty }e^{-\theta k}U(k)V(-\varphi
(n)-k)  \notag \\
&&\qquad \leq C_{2}b_{n}V(-\varphi (n))\left( \sum_{k=1}^{\infty }e^{-\theta
k}U(k)(1+V(-2k))\right)  \notag \\
&&\qquad \leq C_{3}b_{n}V(-\varphi (n)).  \label{First}
\end{eqnarray}

We know by (\ref{AsymV}) that $V(-w)\in R_{\alpha (1-\rho )}.$ Therefore,%
\begin{equation}
\int_{0}^{y}V(-w)dw\sim \frac{yV(-y)}{\alpha (1-\rho )+1}\in R_{\alpha
(1-\rho )+1}  \label{Regular3}
\end{equation}%
as $y\rightarrow \infty $ (see \cite[Ch.VIII, Sec. 9, Theorem 1]{Fel71}).
This fact combined with (\ref{First}) and Corollary~\ref{C_IntegVW} shows
that
\begin{eqnarray}
&&\sum_{j=\left[ n/2\right] +1}^{n}\mathbf{E}\left[ e^{\theta
S_{j}};S_{n}\leq \varphi (n),\tau _{n}=j\right] \leq Cb_{n}V(-\varphi (n))
\notag \\
&&\qquad \qquad =o\left( b_{n}\int_{0}^{\varphi (n)}V(-w)dw\right) =o(%
\mathbf{P}\left( S_{n}\leq \varphi (n),L_{n}\geq 0\right) )  \label{Big_j}
\end{eqnarray}%
as $n\rightarrow \infty ,$ and that there exists a constant $C>0$ such that
\begin{equation*}
\int_{0}^{2y}V(-w)dw\leq C\int_{0}^{y}V(-w)dw
\end{equation*}%
for all sufficiently large $y$. Combining this estimate with (\ref{RightTail}%
) we conclude that
\begin{eqnarray}
&&\sum_{J+1\leq j\leq n/2}\mathbf{E}\left[ e^{\theta S_{j}};S_{n}\leq
\varphi (n),\tau _{n}=j\right]  \notag \\
&&\qquad \leq C\sum_{J\leq j\leq n/2}b_{j}b_{n-j}\sum_{k=1}^{\infty
}e^{-\theta k}U(k)\int_{0}^{\varphi (n)+k}V(-w)dw  \notag \\
&&\qquad \leq C_{1}b_{n}\int_{0}^{2\varphi (n)}V(-w)dw\sum_{J\leq j\leq
n/2}b_{j}\left( \sum_{k=1}^{\infty }e^{-\theta k}U(k)\left(
1+\int_{0}^{2k}V(-w)dw\right) \right)  \notag \\
&&\qquad \leq C_{2}\sum_{J\leq j\leq \infty }b_{j}\times
b_{n}\int_{0}^{\varphi (n)}V(-w)dw=\varepsilon _{J}\mathbf{P}\left(
S_{n}\leq \varphi (n),L_{n}\geq 0\right) ,  \label{Small_j}
\end{eqnarray}%
where $\varepsilon _{J}\rightarrow 0$ as $J\rightarrow \infty $ in view of (%
\ref{Defb}) and Corollary \ref{C_IntegVW}. Estimates (\ref{Big_j}) and (\ref%
{Small_j}) imply
\begin{equation}
\lim_{J\rightarrow \infty }\lim_{n\rightarrow \infty }\frac{R(J+1,n)}{%
\mathbf{P}\left( S_{n}\leq \varphi (n),L_{n}\geq 0\right) }\leq
C\lim_{J\rightarrow \infty }\varepsilon _{J}=0.  \label{Negl1}
\end{equation}

We now take $j\in \lbrack 0,J]$ and use the decomposition%
\begin{eqnarray*}
\mathbf{E}\left[ e^{\theta S_{j}};S_{n}\leq \varphi (n),\tau _{n}=j\right]
&=&\mathbf{E}\left[ e^{\theta S_{j}};S_{n}\leq \varphi (n),\tau
_{n}=j,S_{j}>-\sqrt{\varphi (n)}\right] \\
&&+\mathbf{E}\left[ e^{\theta S_{j}};S_{n}\leq \varphi (n),\tau
_{n}=j,S_{j}\leq -\sqrt{\varphi (n)}\right] .
\end{eqnarray*}%
According to Lemma 5 in \cite{VD2023}, for each fixed $j$
\begin{equation*}
\mathbf{E}\left[ e^{\theta S_{j}};S_{n}\leq \varphi (n),\tau
_{n}=j,S_{j}\leq -\sqrt{\varphi (n)}\right] =o(\mathbf{P}\left( S_{n-j}\leq
\varphi (n),L_{n-j}\geq 0\right) )
\end{equation*}%
as $n\rightarrow \infty $. Further,
\begin{eqnarray*}
&&\mathbf{E}\left[ e^{\theta S_{j}};S_{n}\leq \varphi (n),\tau _{n}=j,S_{j}>-%
\sqrt{\varphi (n)}\right] \\
&&\qquad =\int_{-\sqrt{\varphi (n)}}^{0}e^{\theta x}\mathbf{P}\left(
S_{j}\in dx,\tau _{j}=j\right) \mathbf{P}(S_{n-j}\leq \varphi
(n)-x,L_{n-j}\geq 0).
\end{eqnarray*}%
For $x\in \lbrack -\sqrt{\varphi (n)},0]$ and fixed $j$ we have
\begin{equation*}
0\leq \varphi (n)-x\leq \varphi (n)+\sqrt{\varphi (n)}=o(a_{n})=o(a_{n-j})
\end{equation*}%
as $n\rightarrow \infty $. This estimate and Corollary \ref{C_IntegVW} imply
\begin{equation*}
\mathbf{P}(S_{n-j}\leq \varphi (n)-x,L_{n-j}\geq 0)\sim g_{\alpha ,\beta
}(0)b_{n-j}\int_{0}^{\varphi (n)-x}V(-w)dw
\end{equation*}
as $n\rightarrow \infty $ uniformly in $x\in (-\sqrt{\varphi (n)},0]$. We
know that $b_{n}\in R_{1+\alpha ^{-1}}$. Therefore, $b_{n-j}\sim b_{n}$ as $%
n\rightarrow \infty $ for each fixed $j$. Besides,
\begin{equation*}
\lim_{n\rightarrow \infty }\sup_{x\in (-\sqrt{\varphi (n)},0]}\left( \frac{%
\int_{0}^{\varphi (n)-x}V(-w)dw}{\int_{0}^{\varphi (n)}V(-w)dw}-1\right) =0
\end{equation*}%
according to (\ref{Regular3}). Hence we conclude that, for each fixed $j$
\begin{eqnarray*}
\mathbf{P}(S_{n-j} \leq \varphi (n)-x,L_{n-j}\geq 0)&\sim& g_{\alpha ,\beta
}(0)b_{n}\int_{0}^{\varphi (n)}V(-w)dw \\
&\sim &\mathbf{P}(S_{n}\leq \varphi (n),L_{n}\geq 0)
\end{eqnarray*}%
as $n\rightarrow \infty $ uniformly in $x\in (-\sqrt{\varphi (n)},0]$. Thus,
as $n\rightarrow \infty $
\begin{eqnarray*}
&&\int_{-\sqrt{\varphi (n)}}^{0}e^{\theta x}\mathbf{P}\left( S_{j}\in
dx,\tau _{j}=j\right) \mathbf{P}(S_{n-j}\leq \varphi (n)-x,L_{n-j}\geq 0) \\
&&\qquad \qquad \sim \mathbf{P}(S_{n}\leq \varphi (n),L_{n}\geq
0)\int_{-\infty }^{0}e^{\theta x}\mathbf{P}\left( S_{j}\in dx,\tau
_{j}=j\right) .
\end{eqnarray*}%
Recalling (\ref{Negl1}) and summing the estimates above over $j$ from $0$ to
$\infty $ we get by~(\ref{Decomp1}) that
\begin{eqnarray*}
\mathbf{E}\left[ e^{S_{\tau _{n}}};S_{n}\leq \varphi (n)\right] &\sim &%
\mathbf{P}(S_{n}\leq \varphi (n),L_{n}\geq 0)\int_{-\infty }^{0}e^{\theta
x}\sum_{j=0}^{\infty }\mathbf{P}\left( S_{j}\in dx,\tau _{j}=j\right) \\
&=&\mathbf{P}(S_{n}\leq \varphi (n),L_{n}\geq 0)\left( 1+\int_{-\infty
}^{0}e^{\theta x}U(-dx)\right) \\
&=&\mathbf{P}(S_{n}\leq \varphi (n),L_{n}\geq 0)\theta \int_{0}^{\infty
}e^{-\theta x}U(x)dx
\end{eqnarray*}%
as $n\rightarrow \infty $, where at the last step we have used the equality $%
U(0+)=1$.

We complete the proof of Theorem \ref{T_minimPosit} by combining the
obtained equivalence relation with Corollary \ref{C_IntegVW}.

\section{Proof of Theorem \ref{T_lowend}}

\label{Sec4}

We write
\begin{eqnarray*}
&&\mathbf{E}_{x}\left[ e^{\theta S_{n}};S_{n}\leq \psi (n),M_{n}<0\right]
=\int_{-\infty }^{\psi (n)}e^{\theta w}\mathbf{P}_{x}\left( S_{n}\in
dw,M_{n}<0\right) \\
&&\qquad \qquad \qquad \quad =e^{\theta \psi (n)}\int_{0}^{\infty
}e^{-\theta y}\mathbf{P}_{x}\left( S_{n}\in \psi (n)-dy,M_{n}<0\right) .
\end{eqnarray*}%
For any $h\in (0,1]$ and $N\in \mathbb{N}$ we have%
\begin{eqnarray*}
e^{-\theta h}T_{1}(n,N,h) &\leq &\int_{0}^{\infty }e^{-\theta y}\mathbf{P}%
_{x}\left( S_{n}\in \psi (n)-dy,M_{n}<0\right) \\
&\leq &T_{1}(n,N,h)+T_{2}(n,N),
\end{eqnarray*}%
where%
\begin{eqnarray*}
T_{1}(n,N,h):= &&\sum_{k=0}^{\left[ N/h\right] }e^{-\theta kh}\mathbf{P}%
_{x}\left( S_{n}\in \lbrack \psi (n)-\left( k+1\right) h,\psi
(n)-kh),M_{n}<0\right) , \\
T_{2}(n,N):= &&\sum_{r=N}^{\infty }e^{-\theta r}\mathbf{P}_{x}\left(
S_{n}\in \lbrack \psi (n)-\left( r+1\right) ,\psi (n)-r),M_{n}<0\right) .
\end{eqnarray*}%
Since $U(w),\,w\geq 0,$ is a renewal function, there exists a constant $%
C_{1}\geq 1$ such that $U(w)\leq C_{1}(w+1)$ for all $w\geq 0$.

Combining this inequality with (\ref{LocalMaxim}) shows that there exists $%
N_{0}$ such that the estimates
\begin{eqnarray*}
T_{2}(n,N) &\leq &C\,b_{n}V(x)\sum_{r=N}^{\infty }e^{-\theta r}\,U(r+1-\psi
(n)) \\
&\leq &C\,b_{n}V(x)U(-\psi (n))\left( \sum_{r=N}^{\infty }e^{-\theta
r}\,(1+U(2r+2))\right) \\
&=&C\,b_{n}V(x)U(-\psi (n))e^{-\theta N}\left( \sum_{r=0}^{\infty
}e^{-\theta r}\,(1+U(2r+2N+2))\right) \\
&\leq &2C\,b_{n}V(x)U(-\psi (n))e^{-\theta N}\sum_{r=0}^{\infty }e^{-\theta
r}\,(r+N+2)\, \\
&\leq &C\,_{1}b_{n}V(x)U(-\psi (n))Ne^{-\theta N}
\end{eqnarray*}%
are valid for all $N\geq N_{0}$.

Further, if Condition A1 is satisfied, then, according to Corollary \ref%
{C_Doney_regular}
\begin{eqnarray}
&&\mathbf{P}_{x}\left( S_{n}\in \lbrack \psi (n)-\left( k+1\right) h,\psi
(n)-kh),M_{n}<0\right)  \notag \\
&&\qquad \qquad \qquad \sim g_{\alpha ,-\beta }(0)b_{n}V(x)hU((k+1)h-\psi
(n))  \label{DDon}
\end{eqnarray}%
as $n\rightarrow \infty $ uniformly in negative $x=o(a_{n})$ and $\psi
(n)-\left( k+1\right) h=o(a_{n})$ with $\left( k+1\right) h\leq 2N$.
Besides,
\begin{equation*}
\frac{U(z-\psi (n))}{U(-\psi (n))}\rightarrow 1
\end{equation*}%
as $n\rightarrow \infty $ uniformly in $z=o(\psi (n))$. As a result, we get
\begin{eqnarray*}
T_{1}(n,N,h) &=&(1+o(1))g_{\alpha ,-\beta }(0)b_{n}V(x)h\sum_{k=0}^{\left[
N/h\right] }e^{-\theta kh}U(kh-\psi (n)) \\
&=&(1+o(1))g_{\alpha ,-\beta }(0)b_{n}V(x)U\left( -\psi (n)\right)
h\sum_{k=0}^{\left[ N/h\right] }e^{-\theta kh} \\
&\leq &(1+\varepsilon )g_{\alpha ,-\beta }(0)b_{n}V(x)U\left( -\psi
(n)\right) h(1-e^{-\theta h})^{-1}
\end{eqnarray*}%
for any $\varepsilon >0$ and sufficiently large $n\geq n_{0}(\varepsilon )$.
Thus,
\begin{equation*}
\limsup_{n\rightarrow \infty }\frac{\int_{0}^{\infty }e^{-\theta y}\mathbf{P}%
_{x}\left( S_{n}\in \psi (n)-dy,M_{n}<0\right) }{b_{n}U(-\psi (n))V(x)}\leq
(1+\varepsilon )g_{\alpha ,-\beta }(0)h(1-e^{\theta h})^{-1}.
\end{equation*}%
Since $\varepsilon >0$ and $h>0$ may be selected arbitrary small, we
conclude that%
\begin{equation*}
\limsup_{n\rightarrow \infty }\frac{\int_{0}^{\infty }e^{-\theta y}\mathbf{P}%
_{x}\left( S_{n}\in \psi (n)-dy,M_{n}<0\right) }{b_{n}U(-\psi (n))V(x)}\leq
\frac{g_{\alpha ,-\beta }(0)}{\theta }.
\end{equation*}%
In a similar way we get for any $\varepsilon >0$ and sufficienly large $n$%
\begin{eqnarray*}
e^{-\theta h}T_{1}(n,N,h) &\geq &(1-\varepsilon )g_{\alpha ,-\beta
}(0)b_{n}V(x)U\left( -\psi (n)\right) h\times \\
&&\times \left( (1-e^{-\theta h})^{-1}-\sum_{r=N}^{\infty }e^{-\theta
r}\right) ,
\end{eqnarray*}%
implying that%
\begin{eqnarray*}
&&\liminf_{n\rightarrow \infty }\frac{\int_{0}^{\infty }e^{-\theta y}\mathbf{%
P}_{x}\left( S_{n}\in \psi (n)-dy,M_{n}<0\right) }{b_{n}U(-\psi (n))V(x)} \\
&&\qquad \qquad \geq \liminf_{N\rightarrow \infty }\liminf_{h\downarrow
0}\lim_{\varepsilon \downarrow 0}\liminf_{n\rightarrow \infty }\frac{%
e^{-\theta h}T_{1}(n,N,h)}{b_{n}U(-\psi (n))V(x)}\geq \frac{g_{\alpha
,-\beta }(0)}{\theta }.
\end{eqnarray*}%
It follows that%
\begin{equation*}
\lim_{n\rightarrow \infty }\frac{\int_{0}^{\infty }e^{-\theta y}\mathbf{P}%
_{x}\left( S_{n}\in \psi (n)-dy,M_{n}<0\right) }{b_{n}U(-\psi (n))V(x)}=%
\frac{g_{\alpha ,-\beta }(0)}{\theta },
\end{equation*}%
as desired.

\section{Proof of Theorem \ref{T_minim}}

\label{Sec5}

We write%
\begin{equation*}
\mathbf{E}\left[ e^{\theta S_{\tau _{n}}};S_{n}\leq \psi (n)\right]
=R(0,n-J)+\sum_{j=n-J+1}^{n}\mathbf{E}\left[ e^{\theta S_{j}};S_{n}\leq \psi
(n),\tau _{n}=j\right] ,
\end{equation*}%
where%
\begin{equation*}
R(0,n-J):=\sum_{j=0}^{n-J}\mathbf{E}\left[ e^{\theta S_{j}};S_{n}\leq \psi
(n),\tau _{n}=j\right] .
\end{equation*}%
Clearly, for each $j\in \lbrack 0,n]$%
\begin{eqnarray}
&&\mathbf{E}\left[ e^{\theta S_{j}};S_{n}\leq \psi (n),\tau _{n}=j\right]
\notag \\
&&\quad =\int_{-\infty }^{\psi (n)}e^{\theta x}\mathbf{P}\left( S_{j}\in
dx,\tau _{j}=j\right) \mathbf{P}\left( S_{n-j}\leq \psi (n)-x,L_{n-j}\geq
0\right)  \notag \\
&&\qquad =e^{\theta \psi (n)}T_{3}(n,j),  \label{DefT_3}
\end{eqnarray}%
where%
\begin{equation*}
T_{3}(n,j):=\int_{0}^{\infty }e^{-\theta y}\mathbf{P}\left( S_{j}\in \psi
(n)-dy,\tau _{j}=j\right) \mathbf{P}\left( S_{n-j}\leq y,L_{n-j}\geq
0\right) .
\end{equation*}%
In view of (\ref{AsymV}) there is a constant $C\in (0,\infty )$ such that
\begin{equation}
e^{-\theta y}\int_{0}^{y}V(-z)dz\leq Ce^{-\theta y/2}  \label{EstimIntegr}
\end{equation}%
for all $y\geq 0$. Hence, using the inequalities (\ref{Rough1}), (\ref%
{LocalMaxim}) and (\ref{AsymU}) we deduce the estimates
\begin{eqnarray}
&&T_{3}(n,j)\leq Cb_{n-j}\int_{0}^{\infty }e^{-\theta y}\mathbf{P}\left(
S_{j}\in \psi (n)-dy,\tau _{j}=j\right) \int_{0}^{y}V(-z)dz  \notag \\
&&\, \leq Cb_{n-j}\sum_{k=0}^{\infty }e^{-\theta k}\mathbf{P}\left(
S_{j}\in \lbrack \psi (n)-k-1,\psi (n)-k),\tau _{j}=j\right)
\int_{0}^{k+1}V(-z)dz  \notag \\
&&\, \leq C_{1}b_{n-j}b_{j}\sum_{k=0}^{\infty }e^{-\theta k/2}U(k+1-\psi
(n))  \notag \\
&&\, \leq C_{1}b_{n-j}b_{j}U(-2\psi (n))\sum_{k=0}^{\infty }e^{-\theta
k/2}\left( U(2k+2)+1\right)  \notag \\
&&\, \leq C_{2}b_{n-j}b_{j}U(-\psi (n)).  \label{MediumTerm}
\end{eqnarray}%
In view of (\ref{Defb}) $b_{n-j}\leq Cb_{n}$ for all $j\leq n/2$. Therefore,%
\begin{eqnarray*}
&&\sum_{j=J}^{n-J}\mathbf{E}\left[ e^{\theta S_{j}};S_{n}\leq \psi (n),\tau
_{n}=j\right] \leq C_{2}U(-\psi (n))e^{\theta \psi
(n)}\sum_{j=J}^{n-J}b_{n-j}b_{j} \\
&&\qquad \quad \leq C_{3}b_{n}U(-\psi (n))e^{\theta \psi
(n)}\sum_{j=J}^{n/2}b_{j}\leq \varepsilon _{J}b_{n}U(-\psi (n))e^{\theta
\psi (n)},
\end{eqnarray*}%
where $\varepsilon _{J}=C_{3}\sum_{j=J}^{\infty }b_{j}\rightarrow 0$ as $%
J\rightarrow \infty $. Further, for any fixed $j\leq J$ we have%
\begin{eqnarray*}
T_{3}(n,j) &\leq &Cb_{n-j}\sum_{k=0}^{\infty }e^{-\theta k/2}\mathbf{P}%
\left( S_{j}\in \lbrack \psi (n)-k-1,\psi (n)-k),\tau _{j}=j\right) \\
&\leq &C_{1}b_{n}\mathbf{P}\left( S_{j}\leq \psi (n)\right)
\sum_{k=0}^{\infty }e^{-\theta k/2}=o(b_{n})
\end{eqnarray*}%
as $n\rightarrow \infty $. Thus,%
\begin{equation*}
\limsup_{J\rightarrow \infty }\limsup_{n\rightarrow \infty }\frac{R(0,n-J)}{%
b_{n}U\left( -\psi (n)\right) e^{\theta \psi (n)}}=0.
\end{equation*}%
Consider now $j=n-t$ for $t\in \lbrack 0,J].$ Then, for any $N\in \mathbb{N}$
and $h>0$ we have%
\begin{equation*}
T_{3}(n,j)\leq T_{4}(n,N,h,j)+T_{5}(n,N,j),
\end{equation*}%
where%
\begin{equation*}
T_{4}(n,N,h,j):=\sum_{k=0}^{\left[ N/h\right]}e^{-\theta kh}\mathbf{P}%
\left( S_{j}\in\lbrack\psi(n)-(k+1)h,\psi(n)-kh),\tau
_{j}=j\right)\mathbf{P}\left(S_{t}\leq(k+1)h,L_{t}\geq
0\right)
\end{equation*}%
and%
\begin{equation*}
T_{5}(n,N,j):=\sum_{r=N}^{\infty }e^{-\theta r}\mathbf{P}\left( S_{j}\in
\lbrack \psi (n)-r-1,\psi (n)-r),\tau _{j}=j\right) \mathbf{P}\left(
S_{t}\leq r+1,L_{t}\geq 0\right) .
\end{equation*}%
Using (\ref{LocalMaxim}) with $x=0$ we conclude that
\begin{eqnarray}
T_{5}(n,N,j) &\leq &C\,b_{j}\sum_{r=N}^{\infty }e^{-\theta r}\,U(r+1-\psi
(n))\mathbf{P}\left( S_{t}\leq r+1,L_{t}\geq 0\right)  \notag \\
&\leq &C\,b_{j}U(-2\psi (n))\sum_{r=N}^{\infty }e^{-\theta r}\,\left(
U(2r+2)+1\right)  \notag \\
&\leq &C_{1}\,b_{j}U(-\psi (n))e^{-\theta N}\sum_{r=0}^{\infty }e^{-\theta
r}\,\left( U(2r+2N+2)+1\right)  \notag \\
&\leq &C_{2}\,b_{j}U(-\psi (n))e^{-\theta N}\sum_{r=0}^{\infty }e^{-\theta
r}\,(r+N+2)  \notag \\
&\leq &C_{3}\,b_{j}U(-\psi (n))Ne^{-\theta N}.  \label{T5_above}
\end{eqnarray}

Further, if $0\leq n-j\leq J$ and $n\rightarrow \infty ,$ then, by (\ref%
{DDon1}) and the duality principle for random walks
\begin{eqnarray*}
&&\mathbf{P}\left( S_{j}\in \lbrack \psi (n)-\left( k+1\right) h,\psi
(n)-kh),\tau _{j}=j\right) \\
&&\qquad =\mathbf{P}\left( S_{j}\in \lbrack \psi (n)-\left( k+1\right)
h,\psi (n)-kh),M_{j}<0\right) \\
&&\qquad \sim g_{\alpha ,-\beta }(0)b_{j}\int_{kh-\psi (n)}^{\left(
k+1\right) h-\psi (n)}U(z)dz\sim g_{\alpha ,-\beta }(0)b_{j}hU(-\psi (n))
\end{eqnarray*}%
as $n\rightarrow \infty $ uniformly in $0\leq \left( k+1\right) h\leq N$.
Using these estimates we see that, for each fixed $h>0$
\begin{equation}
T_{4}(n,N,h,j)\sim g_{\alpha ,-\beta }(0)b_{j}U(-\psi (n))h\sum_{k=0}^{\left[
N/h\right] }e^{-\theta kh}\mathbf{P}\left( S_{t}\leq \left( k+1\right)
h,L_{t}\geq 0\right) .  \label{T4_equiv}
\end{equation}%
Note now that
\begin{eqnarray}
&&h\sum_{k=0}^{\left[ N/h\right] }e^{-\theta kh}\mathbf{P}\left( S_{t}\leq
\left( k+1\right) h,L_{t}\geq 0\right)\notag\\
 &&\qquad\qquad\leq e^{2\theta h}\sum_{k=0}^{\left[
N/h\right] }\int_{\left( k+1\right) h}^{\left( k+2\right) h}e^{-\theta z}%
\mathbf{P}\left( S_{t}\leq z,L_{t}\geq 0\right) dz  \notag \\
&&\qquad\qquad\leq e^{2\theta h}\int_{0}^{\infty }e^{-\theta z}\mathbf{P}\left(
S_{t}\leq z,L_{t}\geq 0\right) dz  \label{H_above}
\end{eqnarray}%
and%
\begin{eqnarray}
&&h\sum_{k=0}^{\left[ N/h\right] }e^{-\theta kh}\mathbf{P}\left( S_{t}\leq
\left( k+1\right) h,L_{t}\geq 0\right)\notag\\
&&\qquad \qquad \geq e^{-2\theta h}\sum_{k=0}^{\left[
N/h\right] }he^{-\theta (k-1)h}\mathbf{P}\left( S_{t}\leq kh,L_{t}\geq
0\right)  \notag \\
&&\qquad \qquad \geq e^{-2\theta
h}\int_{0}^{N-1}e^{-\theta z}\mathbf{P}\left( S_{t}\leq z,L_{t}\geq 0\right)
dz.  \label{H_below}
\end{eqnarray}

Combining (\ref{T4_equiv})-(\ref{H_below}) we see that, for $j\in (n-J,n]$
\begin{equation*}
\lim_{h\downarrow 0}\lim_{N\rightarrow \infty }\lim_{n\rightarrow \infty }%
\frac{T_{4}(n,N,h,j)}{g_{\alpha ,-\beta }(0)b_{j}U(-\psi (n))}%
=\int_{0}^{\infty }e^{-\theta z}\mathbf{P}\left( S_{t}\leq z,L_{t}\geq
0\right) dz.
\end{equation*}%
Since $b_{j}\sim b_{n}$ for $j\in (n-J,n]$, we deduce on account of (\ref%
{T5_above}) that%
\begin{equation}
\limsup_{n\rightarrow \infty }\frac{T_{3}(n,j)}{g_{\alpha ,-\beta
}(0)b_{n}U(-\psi (n))}\leq \int_{0}^{\infty }e^{-\theta z}\mathbf{P}\left(
S_{t}\leq z,L_{t}\geq 0\right) dz.  \label{Upper1}
\end{equation}%
To get the same estimate from below we observe that%
\begin{eqnarray*}
&&T_{3}(n,j)\geq T_{6}(N,n,h,j)\notag\\
&&:=\sum_{k=1}^{\left[ N/h\right] }e^{-\theta (k+1)h}\mathbf{P%
}\left( S_{j}\in \lbrack \psi (n)-\left( k+1\right) h,\psi (n)-kh),\tau
_{j}=j\right) \mathbf{P}\left( S_{t}\leq kh,L_{t}\geq 0\right) \\
&&\sim g_{\alpha ,-\beta }(0)b_{j}U(-\psi (n))h\sum_{k=1}^{\left[ N/h\right]
}e^{-\theta (k+1)h}\mathbf{P}\left( S_{t}\leq kh,L_{t}\geq 0\right)
\end{eqnarray*}%
as $n\rightarrow \infty $ and
\begin{eqnarray*}
&&h\sum_{k=1}^{\left[N/h\right]}e^{-\theta(k+1)h}\mathbf{P}\left(S_{t}\leq
kh,L_{t}\geq 0\right)\geq e^{-2\theta h}\sum_{k=1}^{\left[ N/h\right]
}he^{-\theta (k-1)h}\mathbf{P}\left( S_{t}\leq kh,L_{t}\geq 0\right) \\
&&\qquad\qquad\qquad\qquad\qquad\qquad\qquad\qquad\geq e^{-2\theta h}\int_{0}^{N-1}e^{-\theta z}\mathbf{P}\left( S_{t}\leq
z,L_{t}\geq 0\right) dz.
\end{eqnarray*}%
Hence, similarly to (\ref{Upper1}) we get%
\begin{eqnarray*}
\liminf_{n\rightarrow \infty }\frac{T_{3}(n,j)}{g_{\alpha ,-\beta
}(0)b_{j}U(-\psi (n))} &\geq &\lim_{h\downarrow 0}\lim_{N\rightarrow \infty
}\lim_{n\rightarrow \infty }\frac{T_{6}(n,N,h,j)}{g_{\alpha ,-\beta
}(0)b_{j}U(-\psi (n))} \\
&=&\int_{0}^{\infty }e^{-\theta z}\mathbf{P}\left( S_{t}\leq z,L_{t}\geq
0\right) dz.
\end{eqnarray*}%
Thus, for each $t=n-j$%
\begin{equation*}
\lim_{n\rightarrow \infty }\frac{T_{3}(n,j)}{g_{\alpha ,-\beta
}(0)b_{n}U(-\psi (n))}=\int_{0}^{\infty }e^{-\theta z}\mathbf{P}\left(
S_{t}\leq z,L_{t}\geq 0\right) dz.
\end{equation*}%
Recalling now (\ref{DefT_3}), we conclude that, for each fixed $J$%
\begin{eqnarray*}
&&\sum_{j=n-J+1}^{n}\mathbf{E}\left[ e^{\theta S_{j}};S_{n}\leq \psi
(n),\tau _{n}=j\right] \\
&&\qquad \quad \sim g_{\alpha ,-\beta }(0)b_{n}U(-\psi (n))e^{\theta \psi
(n)}\int_{0}^{\infty }e^{-\theta z}\sum_{t=0}^{J}\mathbf{P}\left( S_{t}\leq
z,L_{t}\geq 0\right) dz.
\end{eqnarray*}%
Letting $J$ tend to infinity we see that
\begin{equation*}
\mathbf{E}\left[ e^{\theta S_{\tau _{n}}};S_{n}\leq \psi (n)\right] \sim
g_{\alpha ,-\beta }(0)b_{n}U\left( -\psi (n)\right) e^{\theta \psi
(n)}\int_{0}^{\infty }e^{-\theta z}V_{0}(-z)dz,
\end{equation*}%
as desired.

\section{Proof of Theorem \ref{T_interm}}

\label{Sec6}

We start by the decomposition%
\begin{equation*}
\mathbf{E}\left[ e^{\theta S_{\tau _{n}}};S_{n}\leq K\right]
=R(0,J)+R(J+1,n-J)+R(n-J+1,n),
\end{equation*}%
where%
\begin{equation*}
R(N_{1},N_{2}):=\sum_{j=N_{1}}^{N_{2}}\mathbf{E}\left[ e^{\theta
S_{j}};S_{n}\leq K,\tau _{n}=j\right] .
\end{equation*}%
By Lemma \ref{L_double} and estimate (\ref{EstimIntegr})
\begin{eqnarray*}
&&\mathbf{E}\left[ e^{\theta S_{j}};S_{n}\leq K,\tau _{n}=j\right] \\
&&\quad =\int_{-\infty }^{K\wedge 0}e^{\theta x}\mathbf{P}\left( S_{j}\in
dx,\tau _{j}=j\right) \mathbf{P}\left( S_{n-j}\leq K-x,L_{n-j}\geq 0\right)
\\
&&\quad =e^{\theta K}\int_{K\vee 0}^{\infty }e^{-\theta y}\mathbf{P}\left(
S_{j}\in K-dy,M_{j}<0\right) \mathbf{P}\left( S_{n-j}\leq y,L_{n-j}\geq
0\right) \\
&\leq &Cb_{n-j}\int_{K\vee 0}^{\infty }e^{-\theta y}\mathbf{P}\left(
S_{j}\in K-dy,M_{j}<0\right) \int_{0}^{y}V(-z)dz \\
&\leq &Cb_{n-j}\sum_{k\geq K\vee 0}e^{-\theta k}\mathbf{P}\left( S_{j}\in
\lbrack K-k-1,K-k),M_{j}<0\right) \int_{0}^{k+1}V(-z)dz \\
&\leq &C_{1}b_{n-j}b_{j}\sum_{k\geq K\vee 0}e^{-\theta k/2}U(k+1-K) \\
&\leq &C_{1}b_{n-j}b_{j}U((-2K)\vee 0)\sum_{k\geq K\vee 0}e^{-\theta
k/2}\left( U(2k+2)+1\right) \leq C_{2}b_{n-j}b_{j}.
\end{eqnarray*}%
Thus,%
\begin{equation*}
R(J+1,n-J)\leq C\sum_{j=J+1}^{n-J}b_{n-j}b_{j}\leq
C_{1}b_{n}\sum_{j=J+1}^{n/2}b_{j}\leq \varepsilon _{J}b_{n},
\end{equation*}%
where $\varepsilon _{J}=C_{1}\sum_{j=J+1}^{\infty }b_{j}\rightarrow 0$ as $%
J\rightarrow \infty $.

Further, for fixed $j\in \lbrack 0,J]$ we write%
\begin{equation*}
\frac{\mathbf{E}\left[ e^{\theta S_{j}};S_{n}\leq K,\tau _{n}=j\right] }{%
b_{n}}=\int_{-\infty }^{K\wedge 0}e^{\theta x}\mathbf{P}\left( S_{j}\in
dx,\tau _{j}=j\right) \frac{\mathbf{P}\left( S_{n-j}\leq K-x,L_{n-j}\geq
0\right) }{b_{n}}.
\end{equation*}%
Note that, for any fixed $x\leq K\wedge 0$
\begin{equation*}
\mathbf{P}(S_{n-j}\leq K-x,L_{n-j}\geq 0)\sim g_{\alpha ,\beta
}(0)b_{n-j}\int_{0}^{K-x}V(-w)dw
\end{equation*}%
as $n\rightarrow \infty $, and, by (\ref{Rough1})\ there is a constant $C\in
(0,\infty )$ such that%
\begin{equation*}
\frac{\mathbf{P}\left( 0\leq S_{n}<K-x,L_{n}\geq 0\right) }{b_{n}}\ \leq \
C\,\,\int_{0}^{K-x}V(-w)dw
\end{equation*}%
for all $n$ and all $x\leq K\wedge 0$. These relations combined with (\ref%
{Regular3}) show that, for any $\theta >0$
\begin{eqnarray*}
&&\int_{-\infty }^{K\wedge 0}e^{\theta x}\mathbf{P}\left( S_{j}\in dx,\tau
_{j}=j\right) \int_{0}^{K-x}V(-w)dw \\
&&\qquad \qquad \leq \int_{-\infty }^{0}e^{\theta x}\mathbf{P}\left(
S_{j}\in dx\right) \int_{0}^{K-x}V(-w)dw \\
&&\qquad \qquad \leq C_{1}\int_{-\infty }^{0}e^{\theta x}\mathbf{P}\left(
S_{j}\in dx\right) \left( \left\vert x\right\vert ^{2}+\left\vert
K\right\vert ^{2}+1\right) <\infty .
\end{eqnarray*}%
Recalling Corollary \ref{C_IntegVW}\ and applying the dominated convergence
theorem we conclude that, for each $j\in \lbrack 0,J]$%
\begin{equation*}
\lim_{n\rightarrow \infty }\frac{\mathbf{E}\left[ e^{\theta S_{j}};S_{n}\leq
K,\tau _{n}=j\right] }{b_{n}}=g_{\alpha ,\beta }(0)\int_{-\infty }^{K\wedge
0}e^{\theta x}\mathbf{P}\left( S_{j}\in dx,\tau _{j}=j\right)
\int_{0}^{K-x}V(-w)dw.
\end{equation*}%
As a result,
\begin{eqnarray}
\lim_{J\rightarrow \infty }\lim_{n\rightarrow \infty }\frac{R(0,J)}{b_{n}}
&=&g_{\alpha ,\beta }(0)\int_{-\infty }^{K\wedge 0}e^{\theta
x}\sum_{j=0}^{\infty }\mathbf{P}\left( S_{j}\in dx,\tau _{j}=j\right)
\int_{0}^{K-x}V(-w)dw  \notag \\
&=&g_{\alpha ,\beta }(0)\int_{-\infty }^{K\wedge 0}e^{\theta
x}U(-dx)\int_{0}^{K-x}V(-w)dw<\infty .  \label{RjFinite}
\end{eqnarray}

To evaluate $R(n-J+1,n)$ we write for $j=n-t\in \lbrack n-J+1,n]$ the
representation%
\begin{equation*}
\frac{\mathbf{E}\left[ e^{\theta S_{j}};S_{n}\leq K,\tau _{n}=j\right] }{%
b_{n}}=\int_{-\infty }^{K\wedge 0}e^{\theta x}\frac{\mathbf{P}\left(
S_{j}\in dx,\tau _{j}=j\right) }{b_{n}}\mathbf{P}\left( S_{t}\leq
K-x,L_{t}\geq 0\right)
\end{equation*}%
and observe that by (\ref{LocalMaxim}) and the duality principle for random
walks the term
\begin{eqnarray}
&&\int_{-\infty }^{-H}e^{\theta x}\frac{\mathbf{P}\left( S_{j}\in dx,\tau
_{j}=j\right) }{b_{n}}\mathbf{P}\left( S_{t}\leq K-x,L_{t}\geq 0\right)
\notag \\
&&\qquad \qquad \leq \mathbf{P}\left( L_{t}\geq 0\right) \int_{-\infty
}^{-H}e^{\theta x}\frac{\mathbf{P}\left( S_{j}\in dx,\tau _{j}=j\right) }{%
b_{n}}  \notag \\
&\leq &\mathbf{P}\left( L_{t}\geq 0\right) \sum_{k=H}^{\infty }e^{-\theta k}%
\frac{\mathbf{P}\left( S_{j}\in \lbrack -k-1,-k),\tau _{j}=j\right) }{b_{n}}
\notag \\
&&\qquad \qquad \leq C_{1}\mathbf{P}\left( L_{t}\geq 0\right)
\sum_{k=H}^{\infty }e^{-\theta k}U(k+1)  \label{Int1}
\end{eqnarray}%
is vanishing as $H\rightarrow \infty $. Thus, we are left with%
\begin{equation*}
\int_{-H}^{K\wedge 0}e^{\theta x}\frac{\mathbf{P}\left( S_{j}\in dx,\tau
_{j}=j\right) }{b_{n}}\mathbf{P}\left( S_{t}\leq K-x,L_{t}\geq 0\right) .
\end{equation*}

To evaluate this term we use Theorem 5.1 in \cite{CC2013} (rewritten in our
notation) according to which
\begin{equation*}
\frac{\mathbf{P}\left( S_{j}\in dx,\tau _{j}=j\right) }{b_{n}dx}\sim
g_{\alpha ,-\beta }(0)U(-x)
\end{equation*}%
as $n\rightarrow \infty $ uniformly in $x<0,\ x=o(a_{n})$ if Condition A2 is
valid. Thus,%
\begin{eqnarray}
&&\lim_{n\rightarrow \infty }\int_{-H}^{K\wedge 0}e^{\theta x}\frac{\mathbf{P%
}\left( S_{j}\in dx,\tau _{j}=j\right) }{b_{n}}\mathbf{P}\left( S_{t}\leq
K-x,L_{t}\geq 0\right)  \notag \\
&&\qquad =g_{\alpha ,-\beta }(0)\int_{-H}^{K\wedge 0}e^{\theta x}U(-x)%
\mathbf{P}\left( S_{t}\leq K-x,L_{t}\geq 0\right) dx.  \label{Int2}
\end{eqnarray}%
Combining (\ref{Int1}) and (\ref{Int2}) we deduce that, for fixed $t=n-j\in
\lbrack 0,J]$%
\begin{equation*}
\lim_{n\rightarrow \infty }\frac{\mathbf{E}\left[ e^{\theta S_{j}};S_{n}\leq
K,\tau _{n}=j\right] }{b_{n}}=g_{\alpha ,-\beta }(0)\int_{-\infty }^{K\wedge
0}e^{\theta x}U(-x)\mathbf{P}\left( S_{t}\leq K-x,L_{t}\geq 0\right) dx.
\end{equation*}%
Summing with respect to $t$ from $0$ to $\infty $ we get%
\begin{equation}
\lim_{J\rightarrow \infty }\lim_{n\rightarrow \infty }\frac{R(n-J+1,n)}{b_{n}%
}=g_{\alpha ,-\beta }(0)\int_{-\infty }^{K\wedge 0}e^{\theta
x}U(-x)V_{0}(K-x)dx<\infty .  \label{RNjFinite}
\end{equation}%
Thus,%
\begin{eqnarray*}
\lim_{n\rightarrow \infty }\frac{\mathbf{E}\left[ e^{\theta S_{\tau
_{n}}};S_{n}\leq K\right] }{b_{n}} &=&g_{\alpha ,\beta }(0)\int_{-\infty
}^{K\wedge 0}e^{\theta x}U(-dx)\int_{0}^{K-x}V(-w)dw \\
&&+g_{\alpha ,-\beta }(0)\int_{-\infty }^{K\wedge 0}e^{\theta
x}U(-x)V_{0}(K-x)dx,
\end{eqnarray*}%
as desired.

\section{Survival probability for branching processes evolving in extremely
unfavorable random environment}

\label{Sec7}

In this section we apply the results obtained for random walks to study the
population size of a critical branching process evolving in unfavorable
random environments. To describe the problems we plan to consider we denote
by $\mathfrak{F}$ $=\left\{ \mathfrak{f}\right\} $ the space of all
probability measures on $\mathbb{N}_{0}:=\{0,1,2,...\}$. For notational
reason, we identify a measure $\mathfrak{f}=\left\{ \mathfrak{f}(\left\{
0\right\} ),\mathfrak{f}(\left\{ 1\right\} ),...\right\} \in $ $\mathfrak{F}$
with the respective probability generating function%
\begin{equation*}
f(s)=\sum_{k=0}^{\infty }\mathfrak{f}(\left\{ k\right\} )s^{k},\quad s\in
\lbrack 0,1],
\end{equation*}%
and make no difference between $\mathfrak{f}$ and $f$. Equipped with the
metric of total variation, $\mathfrak{F}$ $=\left\{ \mathfrak{f}\right\}
=\left\{ f\right\} $ becomes a Polish space. Let
\begin{equation*}
F(s)=\sum_{j=0}^{\infty }F\left( \left\{ j\right\} \right) s^{j},\quad s\in
\lbrack 0,1],
\end{equation*}%
be a random variable taking values in $\mathfrak{F}$, and let
\begin{equation*}
F_{n}(s)=\sum_{j=0}^{\infty }F_{n}\left( \left\{ j\right\} \right)
s^{j},\quad s\in \lbrack 0,1],\quad n\in \mathbb{N},
\end{equation*}%
be a sequence of independent probabilistic copies of $F$. The infinite
sequence $\mathcal{E}=\left\{ F_{n},n\in \mathbb{N}\right\} $ is called a
random environment.

A sequence of nonnegative random variables $\mathcal{Z}=\left\{ Z_{n},\ n\in
\mathbb{N}_{0}\right\} $ specified on\ a probability space $(\Omega ,%
\mathcal{F},\mathbf{P})$ is called a branching process in random environment
(BPRE), if $Z_{0}$ is independent of $\mathcal{E}$ and, given $\mathcal{E}$
the process $\mathcal{Z}$ is a Markov chain with
\begin{equation*}
\mathcal{L}\left( Z_{n}|Z_{n-1}=z_{n-1},\mathcal{E}=(f_{1},f_{2},...)\right)
=\mathcal{L}(\xi _{n1}+\ldots +\xi _{ny_{n-1}})
\end{equation*}%
for all $n\in \mathbb{N}$, $z_{n-1}\in \mathbb{N}_{0}$ and $%
f_{1},f_{2},...\in \mathfrak{F}$, where $\xi _{n1},\xi _{n2},\ldots $ is a
sequence of i.i.d. random variables with distribution $f_{n}.$ Thus, $%
Z_{n-1} $ is the $(n-1)$th generation size of the population of the
branching process and $f_{n}$ is the offspring distribution of an individual
at generation $n-1$.

The sequence
\begin{equation*}
S_{0}=0,\quad S_{n}=X_{1}+...+X_{n},\ n\geq 1,
\end{equation*}%
where $X_{i}=\log F_{i}^{\prime }(1),i=1,2,...,$ is called the associated
random walk for the process $\mathcal{Z}$.

We assume below that $Z_{0}=1$ and impose the following restrictions on the
properties of the BPRE.

\paragraph{Condition B1.}

\emph{The elements of the associated random walk satisfy Conditions A1 and
A2. }

According to the classification of \ BPRE's (see, for instance, \cite{agkv}
and \cite{KV2017}), Condition B1 means that we consider the critical BPRE's.

Our second assumption on the environment concerns reproduction laws of
particles. Set%
\begin{equation*}
\gamma (b)=\frac{\sum_{k=b}^{\infty }k^{2}F\left( \left\{ k\right\} \right)
}{\left( \sum_{i=0}^{\infty }iF\left( \left\{ i\right\} \right) \right) ^{2}}%
.
\end{equation*}

\paragraph{Condition B2.}

\emph{There exist $\varepsilon >0$ and $b\in \,$}$\mathbb{N}$ \emph{such
that\ } \emph{\ }
\begin{equation*}
\mathbf{E}[(\log ^{+}\gamma (b))^{\alpha +\varepsilon }]\ <\ \infty ,
\end{equation*}%
\emph{where }$\log ^{+}x=\log (x\vee 1)$\emph{.}

It is known (see \cite[Theorem 1.1 and Corollary 1.2]{agkv}) that if
Conditions B1 and B2 are valid, then there exist a number $\theta \in
(0,\infty )$ and a sequence $l(1),l(2),...,$ slowly varying at infinity such
that, as $n\rightarrow \infty $%
\begin{equation}
\mathbf{P}\left( Z_{n}>0\right) \sim \theta n^{-(1-\rho )}l(n),
\label{AsymBasic}
\end{equation}%
and for any $x\geq 0$
\begin{eqnarray}
\mathbf{P}\left( Z_{n}>0,S_{n}\leq xa_{n}\right) &=&\mathbf{P}\left(
S_{n}\leq xa_{n}|Z_{n}>0\right) \mathbf{P}\left( Z_{n}>0\right)  \notag \\
&\sim &\mathbf{P}\left( Y_{1}^{+}\leq x\right) \mathbf{P}\left(
Z_{n}>0\right) ,  \label{Meander00}
\end{eqnarray}%
where $\mathcal{Y}^{+}=\left\{ Y_{t}^{+},0\leq t\leq 1\right\} $ denotes the
meander of the strictly stable process $\mathcal{Y}$ with index $\alpha $.

Thus, if a BPRE is critical, then, given $Z_{n}>0$ the random variable $%
S_{n},$ the value of the associated random walk that provides survival of
the population to a distant moment $n$, \ grows like $a_{n}$ times a random
positive multiplier. Since $\mathbf{P}\left( Y_{1}^{+}\leq 0\right) =0$ it
follows that if $\varphi (n)$ satisfies the restriction
\begin{equation}
\limsup_{n\rightarrow \infty }\frac{\varphi (n)}{a_{n}}\leq 0,
\label{Nonfav}
\end{equation}%
then%
\begin{equation*}
\mathbf{P}\left( Z_{n}>0,S_{n}\leq \varphi (n)\right) =\mathbf{P}\left(
S_{n}\leq \varphi (n)|Z_{n}>0\right) \mathbf{P}\left( Z_{n}>0\right)
=o\left( \mathbf{P}\left( Z_{n}>0\right) \right)
\end{equation*}%
as $n\rightarrow \infty $.

It is natural to consider the environment meeting condition (\ref{Nonfav})
as unfavorable for the development of the critical BPRE.

An important case of the unfavorable random environment was considered in
\cite{VD2022}, where was shown that if $\varphi (n)\rightarrow \infty $ as $%
n\rightarrow \infty $ in such a way that $\varphi (n)=o(a_{n})$ then
\begin{equation*}
\mathbf{P}\left( Z_{n}>0,S_{n}\leq \varphi (n)\right) \sim \Theta
b_{n}\int_{0}^{\varphi (n)}V(-w)dw,
\end{equation*}%
where $\Theta \in (0,\infty )$ is a known constant and $V(x)$ is the same as
in (\ref{RenUV}).

In this paper we complement the results of \cite{VD2022} by imposing even
more severe restrictions on the environment. Namely, we assume that the
environment meets the assumption $S_{n}\leq K$ for a fixed constant $K$ and
call such an environment as extremely unfavorable.

Here is the main result of this section.

\begin{theorem}
\label{T_BPRE_interm}Let Conditions B1 and B2 be valid. Then, for any fixed $%
K$
\begin{equation*}
\lim_{n\rightarrow \infty }\frac{\mathbf{P}\left( Z_{n}>0,S_{n}\leq K\right)
}{b_{n}}=G_{left}(K)+G_{right}(K),
\end{equation*}%
where the constants $G_{left}(K)\in (0,\infty )$ and $G_{right}(K)\in
(0,\infty )$ are specified below by formulas (\ref{Gleft}) and (\ref{Gright}%
), respectively.
\end{theorem}

To prove the theorem we introduce two new probability measures $\mathbf{P}%
^{+}$ and $\mathbf{P}^{-}$ using the identities
\begin{equation*}
\begin{array}{rl}
\mathbf{E}[U(x+X);X+x\geq 0]\ =\ U(x)\ , & x\geq 0\ , \\
\mathbf{E}[V(x+X);X+x<0]\ =\ V(x)\ , & x\leq 0\ ,%
\end{array}%
\end{equation*}%
valid for any oscillating random walk (see \cite{KV2017}, Chapter 4.4.3).
The construction procedure of these measures is standard and is explained
for $\mathbf{P}^{+}$ and $\mathbf{P}^{-}$ in detail in \cite{agkv} and \cite%
{ABGV2011} (see also \cite{KV2017}, Chapter 5.2). We recall here only some
basic definitions related to this construction.

Let $\mathcal{F}_{n},n\geq 0,$ be the $\sigma $-field of events generated by
random variables $F_{1},F_{2},...,F_{n}$ and $Z_{0},Z_{1},...,Z_{n}$. These $%
\sigma $-field form a filtration $\mathcal{F}$. We assume that the random
walk $\mathcal{S=}\left\{ S_{n},n\geq 0\right\} $ with the inital value $%
S_{0}=x,x\in \mathbb{R}$, is adapted to filtration $\mathcal{F}$ and
construct probability measures $\mathbf{P}_{x}^{+}$, $x\geq 0,$ as follows.
For every sequence $T_{0},T_{1},...$ of random variables with values in some
space $\mathcal{T}$ and adopted to $\mathcal{F}$ and for any bounded and
measurable function $g:\mathcal{T}^{n+1}\rightarrow \mathbb{R}$, $n\in
\mathbb{N}_{0}$, we set%
\begin{equation*}
\mathbf{E}_{x}^{+}[g(T_{0},\ldots ,T_{n})]\ :=\ \frac{1}{U(x)}\mathbf{E}%
_{x}[g(T_{0},\ldots ,T_{n})U(S_{n});L_{n}\geq 0].\
\end{equation*}%
Similarly, $V$ gives rise to probability measures $\mathbf{P}_{x}^{-}$, $%
x\leq 0$, characterized for each $n\in \mathbb{N}_{0}$ by the equation
\begin{equation*}
\mathbf{E}_{x}^{-}[g(T_{0},\ldots ,T_{n})]\ :=\ \frac{1}{V(x)}\mathbf{E}%
_{x}[g(T_{0},\ldots ,T_{n})V(S_{n});M_{n}<0]\ .
\end{equation*}

The following two lemmas provide the major steps for proving Theorem \ref%
{T_BPRE_interm}.

\begin{lemma}
\label{L_cond} (see Lemma 4 in \cite{VD2022}) Assume Condition $B1$. Let $%
H_{1},H_{2},...,$ be a uniformly bounded sequence of real-valued random
variables adapted to some filtration $\mathcal{\tilde{F}=}\left\{ \mathcal{%
\tilde{F}}_{k},k\in \mathbb{N}\right\} $, which converges $\mathbf{P}^{+}$%
-a.s. to a random variable $H_{\infty }$. Suppose that $\varphi (n),$ $n\in
\mathbb{N},$ is a real-valued function such that $\inf_{n\in \mathbb{N}%
}\varphi (n)\geq C>~0$ and $\varphi (n)=o(a_{n})$ as $n\rightarrow \infty $.
Then%
\begin{equation*}
\lim_{n\rightarrow \infty }\frac{\mathbf{E}\left[ H_{n};S_{n}\leq \varphi
(n),L_{n}\geq 0\right] }{\mathbf{P}\left( S_{n}\leq \varphi (n),L_{n}\geq
0\right) }=\mathbf{E}^{+}\left[ H_{\infty }\right] .
\end{equation*}
\end{lemma}

The statement of the following lemma uses sequences $F_{1},\ldots ,F_{n}$ of
probability generating functions.

\begin{lemma}
\label{L_Subcr} Let $0<\delta <1$. Let
\begin{equation}
W_{n}=g_{n}(F_{1},\ldots ,F_{\lfloor \delta n\rfloor
},Z_{0},Z_{1},...,Z_{\lfloor \delta n\rfloor }),n\geq 1,  \label{Wnew}
\end{equation}%
be random variables with values in an Euclidean (or polish) space $\mathcal{W%
}$ such that
\begin{equation*}
W_{n}\ \rightarrow \ W_{\infty }\quad \mathbf{P}_{x}^{+}\text{-a.s.}
\end{equation*}%
for some $\mathcal{W}$-valued random variable $W_{\infty }$ and for all $%
x\geq 0$. Also let $B_{n}=h_{n}(F_{1},\ldots ,F_{\lfloor \delta n\rfloor
}),n\geq 1$, be random variables with values in an Euclidean (or polish)
space $\mathcal{B}$ such that
\begin{equation*}
B_{n}\ \rightarrow \ B_{\infty }\quad \mathbf{P}^{-}\text{-a.s.}
\end{equation*}%
for some $\mathcal{B}$-valued random variable $B_{\infty }.$ Denote
\begin{equation*}
\tilde{B}_{n}\ :=\ h_{n}(F_{n},\ldots ,F_{n-\lfloor \delta n\rfloor +1})\ .
\end{equation*}%
Then for any bounded, continuous function $\Psi :\mathcal{W}\times \mathcal{B%
}\times \mathbb{R}\rightarrow \mathbb{R}$ and for $\theta >0$ as $%
n\rightarrow \infty $
\begin{align*}
\mathbf{E}[\Psi (W_{n},& \tilde{B}_{n},S_{n})e^{\theta S_{n}}\;;\;\tau
_{n}=n]\;\big/\;\mathbf{E}[e^{\theta S_{n}};\tau _{n}=n] \\
& \rightarrow \ K_{1}\iiint \Psi (w,b,-z)\,\mathbf{P}_{z}^{+}\left(
W_{\infty }\in dw\right) \mathbf{P}^{-}\left( B_{\infty }\in db\right)
U(z)e^{-\theta z}dz\ ,
\end{align*}%
where%
\begin{equation*}
K_{1}^{-1}=K_{1\theta }^{-1}=\int_{0}^{\infty }e^{-\theta z}U(z)\,dz.
\end{equation*}
\end{lemma}

\textbf{Remark.} The statement of the lemma almost literally coincides with
the statement of Theorem 2.8 in \cite{ABGV2011}. The only difference is that
the sequence $\left\{ W_{n},n\geq 1\right\} $ in \cite{ABGV2011} has the
form
\begin{equation*}
W_{n}=g_{n}(F_{1},\ldots ,F_{\lfloor \delta n\rfloor }),n\geq 1.
\end{equation*}%
However, the analysis of the arguments used in \cite{ABGV2011} to check the
validity of Theorem 2.8 shows that they are still valid for $\left\{
W_{n},n\geq 1\right\} $ specified by (\ref{Wnew}).

\textbf{Proof of Theorem \ref{T_BPRE_interm}}. First observe that%
\begin{eqnarray*}
\mathbf{P}\left( Z_{n}>0|\mathcal{E}\right) &\leq &\min_{0\leq k\leq n}%
\mathbf{P}\left( Z_{k}>0|\mathcal{E}\right) \\
&\mathbf{\leq }&\min_{0\leq k\leq n}\mathbf{E}\left[ Z_{k}|\mathcal{E}\right]
=\min_{0\leq k\leq n}e^{S_{k}}=e^{S_{\tau _{n}}}.
\end{eqnarray*}%
This estimate combined with the arguments used in the proof of Theorem \ref%
{T_interm} implies%
\begin{eqnarray}
&&\limsup_{J\rightarrow \infty }\limsup_{n\rightarrow \infty }\frac{\mathbf{P%
}\left( Z_{n}>0,S_{n}\leq K,\tau _{n}\in \left[ J+1,n-J\right] \right) }{%
b_{n}}  \notag \\
&&\quad =\limsup_{J\rightarrow \infty }\limsup_{n\rightarrow \infty }\frac{%
\mathbf{E}\left[ \mathbf{P}\left( Z_{n}>0|\mathcal{E}\right) ;S_{n}\leq
K,\tau _{n}\in \left[ J+1,n-J\right] \right] }{b_{n}}  \notag \\
&&\quad \leq \limsup_{J\rightarrow \infty }\limsup_{n\rightarrow \infty }%
\frac{\mathbf{E}\left[ e^{S_{\tau _{n}}};S_{n}\leq K,\tau _{n}\in \left[
J+1,n-J\right] \right] }{b_{n}}=0.  \label{IntermSmall}
\end{eqnarray}%
Thus, it remains to analyse the case%
\begin{equation*}
\tau _{n}\in \lbrack 0,J]\cup \lbrack n-J+1,n].
\end{equation*}%
We fix sufficiently large positive integers $J$ and $N>\left\vert
K\right\vert $ and, for $j\in \lbrack 0,J]$ consider the chain of estimates%
\begin{eqnarray*}
&&\mathbf{P}\left( Z_{n}>0,S_{n}\leq K,\tau _{n}=j,S_{j}<-N\right) \\
&&\qquad \leq \mathbf{E}\left[ e^{S_{j}};S_{n}\leq K,\tau _{n}=j,S_{j}<-N%
\right] \\
&&\qquad =\int_{-\infty }^{-N}e^{\theta x}\mathbf{P}\left( S_{j}\in dx,\tau
_{j}=j\right) \mathbf{P}\left( S_{n-j}\leq K-x,L_{n-j}\geq 0\right) \\
&&\qquad \leq \int_{-\infty }^{-N}e^{\theta x}\mathbf{P}\left( S_{j}\in
dx,\tau _{j}=j\right) \mathbf{P}\left( S_{n-j}\leq -2x,L_{n-j}\geq 0\right) .
\end{eqnarray*}%
Since $V(-w),w\geq 0,$ is a renewal function, there is a constant $C\in
(0,\infty )$ such that $V(-w)\leq C\left( w+1\right) $. In view of (\ref%
{Rough2}), for any $\varepsilon >0$ one can find sufficiently large $%
N=N(\varepsilon )$ such that
\begin{eqnarray*}
&&\int_{-\infty }^{-N}e^{\theta x}\mathbf{P}\left( S_{j}\in dx,\tau
_{j}=j\right) \mathbf{P}\left( S_{n-j}\leq -2x,L_{n-j}\geq 0\right) \\
&&\qquad \leq Cb_{n-j}\int_{-\infty }^{-N}e^{\theta x}\mathbf{P}\left(
S_{j}\in dx,\tau _{j}=j\right) \int_{0}^{-2x}V(-w)dw \\
&&\qquad \leq C_{1}b_{n-j}\int_{-\infty }^{-N}e^{\theta x}(\left\vert
x\right\vert ^{2}+1)\mathbf{P}\left( S_{j}\in dx\right) \leq \varepsilon
b_{n-j}
\end{eqnarray*}%
for all $n-j\geq 0.$ Thus,%
\begin{equation*}
\mathbf{P}\left( Z_{n}>0,S_{n}\leq K,\tau _{n}=j,S_{j}<-N\right) \leq
\varepsilon b_{n-j}
\end{equation*}%
for sufficiently large $N$.

Further, for any $Q\in \mathbb{N}$ and $N>K$ we get by (\ref{Rough2}) and (%
\ref{AsymV}):
\begin{eqnarray*}
&&\mathbf{P}\left( Z_{n}>0,S_{n}\leq K,\tau _{n}=j,S_{j}\geq
-N,Z_{j}>Q\right) \\
&&\quad \leq \mathbf{P}\left( S_{n}\leq K,\tau _{n}=j,S_{j}\geq
-N,Z_{j}>Q\right) \\
&&\quad =\int_{-N}^{K\wedge 0}\mathbf{P}\left( S_{j}\in dx,\tau
_{j}=j,Z_{j}>Q\right) \mathbf{P}\left( S_{n-j}\leq K-x,L_{n-j}\geq 0\right)
\\
&&\quad \leq \int_{-N}^{0}\mathbf{P}\left( S_{j}\in dx,\tau
_{j}=j,Z_{j}>Q\right) \mathbf{P}\left( S_{n-j}\leq 2N,L_{n-j}\geq 0\right) \\
&&\quad \leq \mathbf{P}\left( \tau _{j}=j,Z_{j}>Q\right) \mathbf{P}\left(
S_{n-j}\leq 2N,L_{n-j}\geq 0\right) \\
&&\quad \leq Cb_{n-j}\mathbf{P}\left( \tau _{j}=j,Z_{j}>Q\right)
\int_{0}^{2N}V(-w)dw \\
&&\quad \leq C_{1}b_{n-j}\mathbf{P}\left( \tau _{j}=j,Z_{j}>Q\right) N^{2}.
\end{eqnarray*}%
Clearly,
\begin{equation*}
\lim_{Q\rightarrow \infty }\mathbf{P}\left( \tau _{j}=j,Z_{j}>Q\right) \leq
\lim_{Q\rightarrow \infty }Q^{-1}\mathbf{E}\left[ e^{S_{j}};\tau _{j}=j%
\right] =0.
\end{equation*}%
Therefore, we can select a sufficiently large $Q$ such that%
\begin{equation*}
\mathbf{P}\left( \tau _{j}=j,Z_{j}>Q\right) N^{2}\leq \varepsilon
\end{equation*}%
and obtain the estimate
\begin{equation*}
\mathbf{P}\left( Z_{n}>0,S_{n}\leq K,\tau _{n}=j,Z_{j}>Q\right) \leq
C_{1}\varepsilon b_{n-j}.
\end{equation*}

We now consider the term%
\begin{equation*}
\mathbf{P}\left( Z_{n}>0,S_{n}\leq K,\tau _{n}=j,S_{j}\geq -N,Z_{j}\leq
Q\right) .
\end{equation*}%
Introduce the event%
\begin{equation*}
A_{u.s}:=\left\{ Z_{n}>0\text{ for all }n\geq 0\right\}
\end{equation*}%
and for $0\leq s\leq 1$ define iterations
\begin{equation*}
F_{k,n}(s):=F_{k+1}(F_{k+2}(\ldots F_{n}(s)\ldots ))\text{ if }0\leq k<n\text{ and }%
F_{n,n}(s):=s.
\end{equation*}%
Since the limit
\begin{equation*}
\lim_{n\rightarrow \infty }\mathbf{P}\left( Z_{n}>0|\mathcal{E}%
,Z_{0}=q\right) =\lim_{n\rightarrow \infty }\left( 1-F_{0,n}^{q}(0)\right)
=:1-F_{0,\infty }^{q}(0)
\end{equation*}%
exists $\mathbf{P}^{+}$-a.s. by monotonicity of the extinction probability
of the branching process, it follows from Lemma \ref{L_cond} that
\begin{equation*}
\mathbf{E}\left[ 1-F_{0,n}^{q}(0)|S_{n}\leq y_{n},L_{n}\geq 0\right]
\rightarrow \mathbf{E}^{+}\left[ 1-F_{0,\infty }^{q}(0)\right] =\mathbf{P}%
_{q}^{+}\left( A_{u.s}\right) >0
\end{equation*}%
as $n\rightarrow \infty $ if $y_{n}=o(a_{n})$. Besides, $\mathbf{P}%
_{q}^{+}\left( A_{u.s}\right) >0$ according to Proposition 3.1 in \cite{agkv}%
.

Now Corollary \ref{C_IntegVW} gives that, for fixed $q\leq Q$ and $x\in
\lbrack -N,K]$%
\begin{eqnarray*}
&&\mathbf{P}\left( Z_{n-j}>0,S_{n-j}\leq K-x,L_{n-j}\geq 0,Z_{0}=q\right) \\
&&\qquad =\mathbf{E}\left[ 1-F_{0,n-j}^{q}(0);S_{n-j}\leq K-x,L_{n-j}\geq 0%
\right] \\
&&\qquad \sim \mathbf{P}_{q}^{+}\left( A_{u.s}\right) \mathbf{P}\left(
S_{n-j}\leq K-x,L_{n-j}\geq 0\right) \\
&&\qquad \sim g_{\alpha ,\beta }(0)\mathbf{P}_{q}^{+}\left( A_{u.s}\right)
b_{n-j}\int_{0}^{K-x}V(-w)dw
\end{eqnarray*}%
as $n-j\rightarrow \infty $. Hence, using estimate (\ref{Rough2}) once again
and applying the dominated convergence theorem we conclude that
\begin{eqnarray*}
&&\frac{\mathbf{P}\left( Z_{n}>0,S_{n}\leq K,\tau _{n}=j,S_{j}\geq
-N,Z_{j}=q\right) }{b_{n-j}} \\
&=&\int_{-N}^{K\wedge 0}\mathbf{P}\left( S_{j}\in dx,\tau
_{j}=j,Z_{j}=q\right) \frac{\mathbf{E}\left[ 1-F_{0,n-j}^{q}(0);S_{n-j}\leq
K-x,L_{n-j}\geq 0\right] }{b_{n-j}} \\
&\sim &g_{\alpha ,\beta }(0)\mathbf{P}_{q}^{+}\left( A_{u.s}\right)
\int_{-N}^{K\wedge 0}\mathbf{P}\left( S_{j}\in dx,\tau _{j}=j,Z_{j}=q\right)
\int_{0}^{K-x}V(-w)dw
\end{eqnarray*}%
as $n-j\rightarrow \infty $.

Combining all these estimates, letting $Q$ and $N$ tend to infinity and
using the equivalence $b_{n-j}\sim b_{n}$ valid for each fixed $j$ as $%
n\rightarrow \infty $, we see that, for any fixed $j\in \lbrack 0,J]$%
\begin{eqnarray*}
m_{j}:= &&\lim_{n\rightarrow \infty }\frac{\mathbf{P}\left(
Z_{n}>0,S_{n}\leq K,\tau _{n}=j\right) }{b_{n}} \\
&=&g_{\alpha ,\beta }(0)\int_{-\infty }^{K\wedge 0}\mathbf{E}\left[ \mathbf{P%
}_{Z_{j}}^{+}\left( A_{u.s}\right) \mathbf{P}\left( S_{j}\in dx,\tau
_{j}=j\right) \int_{0}^{K-x}V(-w)dw\right] .
\end{eqnarray*}%
Thus,%
\begin{equation}
G_{left}(K):=\lim_{J\rightarrow \infty }\lim_{n\rightarrow \infty }\frac{%
\mathbf{P}\left( Z_{n}>0,S_{n}\leq K,\tau _{n}\in \lbrack 0,J]\right) }{b_{n}%
}=\sum_{j=0}^{\infty }m_{j}.  \label{Gleft}
\end{equation}

Clearly,
\begin{equation*}
G_{left}(K)\geq m_{0}=g_{\alpha ,\beta }(0)\mathbf{P}^{+}\left(
A_{u.s}|Z_{0}=1\right) \int_{0}^{K}V(-w)dw>0.
\end{equation*}%
The inequality $G_{left}(K)<\infty $ follows from (\ref{RjFinite}) and the
estimate%
\begin{eqnarray*}
\mathbf{P}\left( Z_{n}>0,S_{n}\leq K,\tau _{n}\in \lbrack 0,J]\right) &=&%
\mathbf{E}\left[ 1-F_{0,n}(0);S_{n}\leq K,\tau _{n}\in \lbrack 0,J]\right] \\
&\leq &\mathbf{E}\left[ e^{S_{\tau _{n}}};S_{n}\leq K,\tau _{n}\in \lbrack
0,J]\right] =R(0,J).
\end{eqnarray*}

We now assume that $n-j=t\in \lbrack 0,J]$ and, using the independency of $%
F_{1},\ldots,F_{j}$ and $F_{j+1},\ldots,F_{n}$, write
\begin{eqnarray*}
&&\mathbf{P}\left( Z_{n}>0,S_{n}\leq K,\tau _{n}=j\right) =\mathbf{E}\left[
1-F_{0,n}(0);S_{n}\leq K,\tau _{n}=j\right] \\
&=&\mathbf{E}\left[ 1-F_{0,j}(F_{j,n}(0));S_{j}\leq K-(S_{n}-S_{j}),\tau
_{n}=j\right] \\
&=&\int_{K}^{\infty }\int_{0}^{1}\mathbf{P}\left( F_{0,t}(0)\in ds,S_{t}\in
dx,L_{t}\geq 0\right) \mathbf{E}\left[ 1-F_{0,j}(s);S_{j}\leq K-x,\tau _{j}=j%
\right] .
\end{eqnarray*}%
Our aim is to investigate the quantity
\begin{equation*}
\mathbf{E}\left[ 1-F_{0,j}(s);S_{j}\leq K-x,\tau _{j}=j\right]
\end{equation*}%
for fixed $s\in \left[ 0,1\right] $ and $x\geq K$.

Clearly,%
\begin{eqnarray*}
\mathbf{E}\left[ 1-s^{Z_{j}}|\mathcal{E}\right] &=&\mathbf{E}\left[ \mathbf{E%
}\left[ 1-s^{Z_{j}}|Z_{[j/2]}\right] \mathcal{E}\right] =\mathbf{E}\left[
1-\left( F_{\left[ j/2\right] ,j}(s)\right) ^{Z_{[j/2]}}|\mathcal{E}\right]
\\
&=&\mathbf{E}\left[ 1-\left( F_{\left[ j/2\right] ,j}(s)\right) ^{\exp
\left\{ S_{\left[ j/2\right] }-S_{j}\right\} Z_{[j/2]}\exp \left\{ -S_{\left[
j/2\right] }\right\} \exp \left\{ S_{j}\right\} }|\mathcal{E}\right] .
\end{eqnarray*}%
Define%
\begin{equation*}
\Psi _{K-x}(w,b,z)=\left( 1-b^{w\exp \left\{ z\right\} }\right)
e^{-z}I\left\{ z\leq K-x\right\}
\end{equation*}%
for $w\geq 0,0\leq b\leq 1,z\in \mathbb{R}$ with $0^{0}=1$ and continue $%
\Psi _{K-x}$ to the other values of $w,b,z$ to a bounded smooth function. In
doing so, points of discontinuity in $(0,0,z)$ and $(w,b,K-x)$ are
unavoidable, which will be bypassed later on.

With this notation in view, we write
\begin{equation*}
\mathbf{E}\left[ 1-F_{0,j}(s);S_{j}\leq K-x,\tau _{j}=j\right] =\mathbf{E}%
\left[ \Psi _{K-x}(W_{j},\tilde{B}_{j}(s),S_{j})e^{S_{j}};\tau _{j}=j\right]
,
\end{equation*}%
where%
\begin{equation*}
W_{j}:=Z_{[j/2]}e^{-S_{\left[ j/2\right] }},\quad \tilde{B}_{j}(s):=\left(
F_{\left[ j/2\right] ,j}(s)\right) ^{\exp \left\{ S_{\left[ j/2\right]
}-S_{j}\right\} }.
\end{equation*}

Since $\left\{ W_{j},j\geq 1\right\} $ is a martingale,%
\begin{equation*}
W_{j}=Z_{[j/2]}e^{-S_{\left[ j/2\right] }}\rightarrow W_{\infty }\qquad
\mathbf{P}_{z}^{+}\text{- a.s.}
\end{equation*}%
as $j\rightarrow \infty $. Moreover, it follows from Proposition 3.1 in \cite%
{agkv} that%
\begin{equation}
\mathbf{P}_{z}^{+}\left( W_{\infty }>0\right) >0  \label{MartinW}
\end{equation}%
for any $z\geq 0$ (In fact, (\ref{MartinW}) was proved in \cite{agkv} for $%
\mathbf{P}_{0}^{+}$ only. However, the analysis of the proof given in the
mentioned proposition shows that (\ref{MartinW}) is valid for any $z>0$ as
well.). Further, according to Lemma 3.2 in \cite{ABGV2011}\ the sequence $%
B_{j}(s):=\left( F_{\left[ j/2\right] ,0}(s)\right) ^{\exp \left\{ -S_{\left[
j/2\right] }\right\} },j=1,2,...$ is nondecreasing and \ $B_{j}(s)\geq
s^{\exp \left\{ -S_{0}\right\} }$. Thus,
\begin{equation*}
B_{j}(s)\rightarrow B_{\infty }(s)\geq s^{\exp \left\{ -S_{0}\right\}
}=s\qquad \mathbf{P}^{-}\text{- a.s.}
\end{equation*}%
as $j\rightarrow \infty $. In view of
\begin{equation*}
\mathbf{P}\left(\tilde{B}_{j}(s)\geq s\right) =\mathbf{P}\left( B_{j}(s)\geq s\right) =1,
\end{equation*}
the second random variable in $\Psi _{K-x}(W_{j},\tilde{B}_{j}(s),S_{j})$ is
separated from zero with probability one.

Thus, we may restrict $\Psi _{K-x}$ to the domain $(w,b,z\neq K-x)$ with $%
b>0,$ where it is continuous. Finally, the additional discontinuity at $%
z=K-x $ has probability $0$ with respect to the measure
\begin{equation*}
\mu \left( dz\right) :=K_{1}e^{-z}U(z)I\left\{ z\geq 0\right\} dz,
\end{equation*}%
where $K_{1}^{-1}:=\int_{0}^{\infty }e^{-z}U(z)dz$.

As a result, we may apply Lemma \ref{L_Subcr} to conclude that, as $%
j\rightarrow \infty $%
\begin{equation}
\frac{\mathbf{E}\left[ \Psi _{K-x}(W_{j},\tilde{B}_{j}(s),S_{j})e^{S_{j}};%
\tau _{j}=j\right] }{\mathbf{E}\left[ e^{S_{j}};\tau _{j}=j\right] }%
\rightarrow h(s,K-x),  \label{Posit2}
\end{equation}%
where%
\begin{equation*}
h(s,K-x):=K_{1}\int \int \int \Psi _{K-x}(w,b,-z)\mathbf{P}_{z}^{+}\left(
W_{\infty }\in dw\right) \mathbf{P}^{-}\left( B_{\infty }(s)\in db\right)
e^{-z}U(z)dz
\end{equation*}%
is a bounded function.

It was shown in the proof of Lemma 3.4 in \cite{ABGV2011} that $B_{\infty
}(s)<1$ $\mathbf{P}^{-}$-a.s. under our conditions. Alltogether this implies
that the right-hand side of (\ref{Posit2}) is positive. Using now the
inequality
\begin{equation*}
\mathbf{E}\left[ 1-F_{0,j}(s);S_{j}\leq K-x,\tau _{j}=j\right] \leq \mathbf{E%
}\left[ e^{S_{j}};\tau _{j}=j\right]
\end{equation*}%
and the dominated convergence theorem we conclude that
\begin{eqnarray}
&&\lim_{j\rightarrow \infty }\int_{K}^{\infty }\int_{0}^{1}\mathbf{P}\left(
F_{0,t}(0)\in ds,S_{t}\in dx,L_{t}\geq 0\right) \frac{\mathbf{E}\left[
1-F_{0,j}(s);S_{j}\leq K-x,\tau _{j}=j\right] }{\mathbf{E}\left[
e^{S_{j}};\tau _{j}=j\right] }  \notag \\
&&\qquad \qquad =\int_{K}^{\infty }\int_{0}^{1}\mathbf{P}\left(
F_{0,t}(0)\in ds,S_{t}\in dx,L_{t}\geq 0\right) h(s,K-x).
\label{SingleTermRight}
\end{eqnarray}%
According to (\ref{MaxSmall})
\begin{equation*}
\mathbf{E}\left[ e^{S_{j}};\tau _{j}=j\right] =\mathbf{E}\left[
e^{S_{j}};M_{j}<0\right] \sim g_{\alpha ,\beta }(0)b_{j}\int_{0}^{\infty
}e^{-w}U(w)dw
\end{equation*}%
as $j\rightarrow \infty $. Using this fact and summing (\ref{SingleTermRight}%
) over $t$ from $0$ to $\infty $ we get
\begin{eqnarray}
&&G_{right}(K):=\lim_{J\rightarrow \infty }\lim_{n\rightarrow \infty }\frac{%
\mathbf{P}\left( Z_{n}>0;S_{n}\leq K,\tau _{n}\in \lbrack n-J+1,n]\right) }{%
b_{n}}  \notag \\
&&\quad =g_{\alpha ,\beta }(0)\int_{0}^{\infty }e^{-w}U(w)dw\times  \notag \\
&&\quad \times \int_{K}^{\infty }\int_{0}^{1}\sum_{t=0}^{\infty }\mathbf{P}%
\left( F_{0,t}(0)\in ds,S_{t}\in dx,L_{t}\geq 0\right) h(s,K-x)>0.
\label{Gright}
\end{eqnarray}%
The inequality $G_{right}(K)<\infty $ follows from (\ref{RNjFinite}) and the
estimate%
\begin{eqnarray*}
&&\mathbf{P}\left( Z_{n}>0,S_{n}\leq K,\tau _{n}\in \lbrack n-J+1,n]\right)
\\
&&\qquad \quad =\mathbf{E}\left[ 1-F_{0,n}(0);S_{n}\leq K,\tau _{n}\in
\lbrack n-J+1,n]\right] \\
&&\qquad \quad \leq \mathbf{E}\left[ e^{S_{\tau _{n}}};S_{n}\leq K,\tau
_{n}\in \lbrack n-J+1,n]\right] =R(n-J+1,n).
\end{eqnarray*}

Combining (\ref{IntermSmall}) with (\ref{Gleft}) and (\ref{Gright}) we get
\begin{equation*}
\lim_{n\rightarrow \infty }\frac{\mathbf{P}\left( Z_{n}>0,S_{n}\leq K\right)
}{b_{n}}=G_{left}(K)+G_{right}(K),
\end{equation*}%
as desired.

\textbf{Acknowledgment. }The work of E.E. Dyakonova and V.A. Vatutin was
performed at the Steklov International Mathematical Center and supported by
the Ministry of Science and Higher Education of the Russian Federation
(agreement no. 075-15-2022-265). The research of C.Dong and V.A. Vatutin was
also supported by the Ministry of Science and Technology of PRC, project
G2022174007L.

\end{document}